\documentclass[12pt]{article}
\usepackage{latexsym}
\usepackage{amsmath}
\usepackage{amssymb}
\usepackage{amsthm}
\usepackage{amscd}
\usepackage[dvipdfmx]{graphicx}
\input cyracc.def
\newfont{\tencyr}{wncyr10}

\begin{document}

\begin{center}
 {\Large \bf   On certain $L$-functions for deformations of knot group representations}
 \end{center}

 \vskip 5pt
 
\begin{center}
Takahiro KITAYAMA, Masanori MORISHITA, Ryoto TANGE\\ and Yuji TERASHIMA
\end{center}

\footnote[0]{
$\;\;\; $ Key words: deformations of representations,  character schemes, twisted knot modules,\\
$\;\;\; $ twisted Alexander invariants, $L$-functions\\
$\;\;\;$ T.K. is partly supported by JSPS Research Fellowships for Young Scientists 26800032\\
$\;\;\; $ M.M. is partly supported by Grants-in-Aid for Scientific Research (B) 24340005.\\
$\;\;\; $ Y.T. is partly supported by Grants-in-Aid for Scientific Research (C)  25400083.
}

\vspace{.5cm}

{\bf Abstract.} We study  the twisted knot module for the universal deformation of an ${\rm SL}_2$-representation
of a knot group, and introduce an associated $L$-function, which may be seen as an analogue of the algebraic $p$-adic $L$-function associated to the Selmer module for the universal deformation of a Galois representation. We then investigate two problems proposed by  Mazur: Firstly we show the torsion property of the twisted knot module over the universal deformation ring under certain conditions. Secondly we compute the $L$-function by some concrete examples for $2$-bridge knots.

\vspace{.8cm}

\begin{center}
{\bf Introduction}
\end{center}

It has been known that there are intriguing analogies between knot theory and number theory (cf. [Mo]).
In particular, it may be  noteworthy that there are close parallels between Alexander-Fox theory and Iwasawa theory ([Ma1], [Mo; Ch.$9\sim 13$]). From the viewpoint of deformations of group representations ([Ma2]), they are concerned with abelian deformations of representations of knot and Galois groups and the associated topological and arithmetic invariants such as the Alexander and Iwasawa polynomials, respectively. In [Ma3], Mazur proposed a number of problems in pursuing these analogies for non-abelian deformations of higher dimensional representations. To carry out Mazur's perspective,  as a first step, we  developed a deformation theory for ${\rm SL}_2$-representations of knot groups in [MTTU]. In this paper, we continue our study and introduce a certain $L$-function
associated to the twisted knot module for the universal deformation of a knot group representation, which may be seen as an analogue of the algebraic $p$-adic $L$-function associated to the Selmer module for the universal deformation of a Galois representation ([G]). 

Let $K$ be a knot in the $3$-sphere $S^3$ and $G_K := \pi_1(S^3 \setminus K)$ the knot group. Fix a field $k$ whose characteristic is not 2 and a complete discrete valuation ring ${\cal O}$ whose residue field is $k$. Let $\overline{\rho} : G_K \rightarrow {\rm SL}_2(k)$ be a given absolutely irreducible representation. It was shown in [MTTU] that there exists the universal deformation ${\boldsymbol \rho} : G_K \rightarrow {\rm SL}_2({\boldsymbol R}_{\overline{\rho}})$ of $\overline{\rho}$, where ${\boldsymbol R}_{\overline{\rho}}$ is a complete local ${\cal O}$-algebra whose residue field is $k$. Assume that ${\boldsymbol R}_{\overline{\rho}}$ is a  Noetherian factorial domain. In this paper, we study the twisted knot module $H_1({\boldsymbol \rho}) := H_1(S^3 \setminus K; {\boldsymbol \rho})$  with coefficients in the universal deformation ${\boldsymbol \rho}$, and introduce the associated  $L$-function $L_K({\boldsymbol \rho})$ defined on the universal deformation space ${\rm Spec}({\boldsymbol R}_{\overline{\rho}})$ as $\Delta_0(H_1({\boldsymbol \rho}))$, the greatest common divisor of generators of the initial Fitting ideal of $H_1({\boldsymbol \rho})$ over the universal deformation ring ${\boldsymbol R}_{\overline{\rho}}$.  In terms of our $H_1({\boldsymbol \rho})$ and $L_K({\boldsymbol \rho})$, we then formulate the problems proposed by Mazur (questions 1 and 2 of [Ma3; page 440] ) as follows.
 \\
$\;\;\;$(1)  Is $H_1({\boldsymbol \rho})$ finitely generated and torsion as an ${\boldsymbol R}_{\overline{\rho}}$-module ?
\\
$\;\;\;$(2) Investigate the order of the zeroes of $L_K({\boldsymbol \rho})$ at prime divisors of 
\\
$\;\;\;\;\;\;\;\;$ ${\rm Spec}({\boldsymbol R}_{\overline{\rho}})$. 

The corresponding problems of (1) and (2) in the arithmetic counterpart, say $(1)_{\rm arith}$ and $(2)_{\rm arith}$ respectively, are important issues in number theory (cf. questions 1 and 2 of [Ma3; page 454]). In fact, $(1)_{\rm arith}$   is a part of the so-called main conjecture for $p$-adic deformations of a Galois representation. For the cyclotomic
 deformation of a Dirichlet character, the affirmative answer to $(1)_{\rm arith}$  is a basic result in Iwasawa theory ([I]), which asserts that the classical Iwasawa module is finitely generated and  torsion over the Iwasawa algebra.
For the Hida deformation (universal ordinary modular ${\rm GL}_2$-deformation) ([H1], [H2]), the affirmative answer to $(1)_{\rm arith}$ has been shown  by Kato and Ochiai ([Kt], [O1], [O2]), which asserts that the dual Selmer module of the Hida deformation is finitely generated and torsion over the the universal ordinary modular deformation ring ($p$-adic Hecke algebra). The problem $(2)_{\rm arith}$  remains  an interesting problem to be explored  and it is related to deep arithmetic issues (cf. question 3 of [Ma3; page 454], and Ribet's theorem on Herbrand's theorem for example [MW]). 

 So it may be interesting  to study the above problems (1) and (2) in the knot theoretic situation. Our results concerning these  are as follows. For (1), we give a criterion for $H_1({\boldsymbol \rho})$ to be finitely generated and torsion over ${\boldsymbol R}_{\overline{\rho}}$ under certain conditions using a  twisted Alexander invariant of $K$ (cf. Theorem 3.2.4). For (2), we  give some concrete examples for 2-bridge knots $K$ such that $L_K({\boldsymbol \rho})$ has only one zero of order 0 or 2 (cf. Example 4.5).
 
 Here are contents of this paper. In Section 1, we recall the deformation theory for ${\rm SL}_2$-representations of a group, which was developed in [MTTU]. In Section 2, we show the relation between the universal deformation ring  and the character scheme over $\mathbb{Z}$ of ${\rm SL}_2$-representations. In Section 3, we study the twisted knot module with coefficients in the universal deformation of  an ${\rm SL}_2$-representation of a knot group, and introduce an associated $L$-function. In Section 4, we discuss some  examples for some 2-bridge knots, for which we study Mazur's problems.
\\
\\
{\it Acknowledgements.} We would like to thank Ted Chinburg, Shinya Harada, Haruzo Hida, Tadashi Ochiai, Jun Ueki and Seidai Yasuda for helpful communications. We also thank the referee for his/her useful comments.\\
\\
{\bf Notation.}  For a local ring $R$, we denote by $\frak{m}_R$ the maximal ideal of $R$. For an integral domain $A$,  we denote by ${\rm char}(A)$ the characteristic of $A$ and by $Q(A)$ the field of fractions of $A$. For $a, b$ in a commutative ring $A$, $a \doteq b$ means $a=bu$ for some unit $u \in A^{\times}$.\\

\begin{center}
{\bf 1. The universal deformation}
\end{center}

In this section, we present a summary of the deformation theory for ${\rm SL}_2$-representations of a group, which was developed in [MTTU]. We also discuss the obstruction to the deformation problem for a  group representation. Throughout this section, let $G$ denote a group.\\
\\
1.1. {\em Pseudo-representations and their deformations.} Let $A$ be  a commutative ring with identity. A map $T : G \rightarrow A$ is called a {\em pseudo-${\rm SL}_2$-representation} over $A$ if the following four conditions are satisfied:
\vspace{.18cm}\\
(P1) $T(e) = 2$  ($e :=$ the identity element of $G$),\\
(P2) $T(g_1g_2) = T(g_2g_1)$ for any $g_1, g_2 \in G$,\\
(P3) $T(g_1)T(g_2)T(g_3) + T(g_1g_2g_3) + T(g_1g_3g_2) -T(g_1g_2)T(g_3)  - T(g_2g_3)T(g_1) - T(g_1g_3)T(g_2)  = 0$ for any $g_1, g_2, g_3\in G$,\\
(P4) $T(g)^2 - T(g^2) = 2$ for any $g \in G$.
\vspace{.18cm}\\
Note that the conditions (P1) $\sim$ (P3) are nothing but Taylor's conditions for a  pseudo-representation of degree 2 ([Ta]) and that (P4) is the condition for determinant 1.  In the following, we say simply a {\em pseudo-representation} for a pseudo-${\rm SL}_2$-representation. The trace ${\rm tr}(\rho)$ of a representation $\rho : G \rightarrow {\rm SL}_2(A)$ satisfies the conditions (P1) $\sim$ (P4) ([Pr, Theorem 4.3]), and, conversely, a pseudo-${\rm SL}_2$-representation is shown to be obtained as the trace of a representation under certain conditions (See Theorem 1.2.1 below). 

Let $k$ be a perfect field and  let $\mathcal{O}$ be a  complete discrete valuation ring with the  residue field $\mathcal{O}/\frak{m}_{\mathcal{O}} = k$.   There is a unique subgroup $V$ of $\mathcal{O}^{\times}$ such that $k^{\times} \simeq V$ and $\mathcal{O}^{\times} = V \times (1+\frak{m}_{\mathcal{O}})$. The composition map $\lambda : k^{\times} \simeq V \hookrightarrow {\mathcal{O}}^{\times}$ is called the {\em Teichm\"{u}ller lift} which satisfies $\lambda(\alpha) \, \mbox{mod}\, \frak{m}_{\mathcal{O}} = \alpha$ for $\alpha \in k$.  It is extended to $\lambda : k \hookrightarrow {\mathcal{O}}$ by $\lambda(0) := 0$. Let $\frak{CL}_{\cal O}$ be the category of complete local $\mathcal{O}$-algebras with residue field $k$. A morphism in $\frak{CL}_{\cal O}$ is an $\mathcal{O}$-algebra homomorphism inducing the identity on residue fields. 

Let $\overline{T} : G \rightarrow k$ be a pseudo-representation over $k$. A couple $(R, T)$ is called an {\em ${\rm SL}_2$-deformation} of $\overline{T}$ if $R \in \frak{CL}_{\cal O}$ and $T : G \rightarrow R$ is a pseudo-${\rm SL}_2$-representation over $R$ such that $T$ mod $\frak{m}_R = \overline{T}$. In the following, we say simply a {\em deformation} of $\overline{T}$ for an ${\rm SL}_2$-deformation. A deformation $({\boldsymbol R}_{\overline{T}}, {\boldsymbol T})$ of $\overline{T}$ is called a {\em universal deformation} if the following universal property is satisfied: ``For any deformation $(R, T)$ of $\overline{T}$ there exists a unique  morphism $\psi :  {\boldsymbol R}_{\overline{T}} \rightarrow  R$ in $\frak{CL}_{\cal O}$ such that $\psi \circ  {\boldsymbol T} = T$."  Namely the correspondence $\psi \mapsto \psi \circ {\boldsymbol T}$ gives the bijection 
$$ {\rm Hom}_{\frak{CL}_{\cal O}}({\boldsymbol R}_{\overline{T}}, R) \simeq \{(R,T) \, | \, \mbox{deformation of} \; \overline{T} \}. $$
By the universal property, a universal deformation $({\boldsymbol R}_{\overline{T}}, {\boldsymbol T})$ of $\overline{T}$ is unique (if it exists) up to isomorphism. The ${\cal O}$-algebra ${\boldsymbol R}_{\overline{T}}$ is called the {\em universal deformation ring} of $\overline{T}$. \\
\\
{\bf Theorem 1.1.1} ([MTTU; Theorem 1.2.1]).  {\em For a pseudo-representation $\overline{T} : G \rightarrow k$, there exists a universal deformation $({\boldsymbol R}_{\overline{T}}, {\boldsymbol T})$ of $\overline{T}$.}\\
\\
We recall the construction of $({\boldsymbol R}_{\overline{T}}, {\boldsymbol T})$.  Let $X_g$ denote a variable indexed by $g \in G$. Then the universal deformation ring ${\boldsymbol R}_{\overline{T}}$ is given by 
$$ {\boldsymbol R}_{\overline{T}} = \mathcal{O}[[X_g \;(g \in G) ]]/{\cal I},$$ 
where $\mathcal{I}$ is the ideal of the formal power series ring $\mathcal{O}[[X_g \;(g \in G) ]]$  generated by the elements of following type:  Setting $T_g := X_g + \lambda(\overline{T}(g))$, 
\vspace{.18cm}\\
(1) $T_e - 2 = X_e +\lambda(\overline{T}(e)) - 2$,\\
(2) $T_{g_1g_2} - T_{g_2g_1} = X_{g_1g_2} - X_{g_2g_1}$,\\
(3) $T_{g_1}T_{g_2}T_{g_3} + T_{g_1g_2g_3} + T_{g_1g_3g_2} -T_{g_1g_2}T_{g_3} - T_{g_2g_3}T_{g_1} - T_{g_1g_3}T_{g_2}$,\\
(4) $T_g^2 - T_{g^2} -2,$
\vspace{.18cm}\\
for $g, g_1, g_2, g_3 \in G$. The universal deformation ${\boldsymbol T} : G \rightarrow {\boldsymbol R}_{\overline{T}}$ is given by 
$${\boldsymbol T}(g) := T_g \;  \mbox{mod} \; \mathcal{I}.$$
Then, for any deformation $(R,T)$ of $\overline{T}$, the morphism $\psi : {\boldsymbol R}_{\overline{T}} \rightarrow R$ in $\frak{CL}_{\cal O}$ defined by $\psi(X_g) := T(g) - \lambda(\overline{T}(g))$ satisfies $\psi \circ {\boldsymbol T} = T$. 

We note that ${\boldsymbol R}_{\overline{T}}$ constructed above is a complete Noetherian local ${\cal O}$-algebra if $G$ is a finitely generated group.\\
\\
1.2. {\em Deformations of an ${\rm SL}_2$-representation.}  We keep the same notations as in 1.1. In this subsection we assume that ${\rm char}(k) \neq 2$, so that $2$ is invertible in $\mathcal{O}$ and hence in any $R \in \frak{CL}_{\cal O}$. Let $\overline{\rho} : G \rightarrow {\rm SL}_2(k)$ be a given representation. We call a couple $(R, \rho)$ an ${\rm SL}_2$-{\em deformation} of $\overline{\rho}$ if $R \in \frak{CL}_{\cal O}$ and $\rho : G \rightarrow {\rm SL}_2(R)$ is a representation such that $\rho$ mod $\frak{m}_R = \overline{\rho}$. In the following, we say simply a {\em deformation} of $\overline{\rho}$ for an ${\rm SL}_2$-deformation. A deformation $({\boldsymbol R}_{\overline{\rho}}, {\boldsymbol \rho})$ of $\overline{\rho}$ is called a {\em universal deformation} of $\overline{\rho}$ if the following universal property is satisfied: ``For any deformation $(R, \rho)$ of $\overline{\rho}$ there exists a unique morphism $\psi : {\boldsymbol R}_{\overline{\rho}} \rightarrow R$ in $\frak{CL}_{\cal O}$ such that $\psi \circ  {\boldsymbol \rho} \approx \rho$". Here two representations $\rho_1, \rho_2 $ of degree 2 over a local ring $A$ are said to be {\em strictly equivalent}, denoted by $\rho_1 \approx \rho_2$, if there is $\gamma \in I_2 + {\rm M}_2(\frak{m}_A)$ such that $\rho_2(g) = \gamma^{-1}\rho_1(g)\gamma$ for all $g \in G$. Namely the correspondence $\psi \mapsto \psi \circ {\boldsymbol \rho}$ gives the bijection 
$$ {\rm Hom}_{\frak{CL}_{\cal O}}({\boldsymbol R}_{\overline{\rho}},R) \simeq \{ (R, \rho)\, | \, \mbox{deformation of}\; \overline{\rho}\}/\approx. $$
By the universal property, a universal deformation $({\boldsymbol R}_{\overline{\rho}}, {\boldsymbol \rho})$ of $\overline{\rho}$ is unique (if it exists) up to strict equivalence. The $\mathcal{O}$-algebra ${\boldsymbol R}_{\overline{\rho}}$ is called the {\em universal deformation ring} of $\overline{\rho}$. 

A deformation  $(R, \rho)$ of $\overline{\rho}$ gives rise to a deformation $(R, {\rm tr}(\rho))$ of the pseudo-representation ${\rm tr}(\overline{\rho}) : G \rightarrow k$. Assume that  $\overline{\rho}$ is absolutely irreducible, namely, the composite of $\overline{\rho}$ with an inclusion ${\rm SL}_2(k) \hookrightarrow {\rm SL}_2(\overline{k})$ is irreducible for an algebraic closure $\overline{k}$ of $k$. Then, by using theorems of Caryaol [Ca; Theorem 1] and Nyssen [Ny; Theorem 1], it can be shown that this correspondence by the trace is indeed bijective. It is here that the condition ${\rm char}(k) \neq 2$ is used. \\
\\
{\bf Theorem 1.2.1} ([MTTU; Theorem 2.1.2]).  {\em Let $\overline{\rho} : G \rightarrow {\rm SL}_2(k)$ be an absolutely irreducible representation and let $R \in \frak{CL}_{\cal O}$. Then the correspondence $\rho \mapsto {\rm tr}(\rho)$ gives the following bijection}:
$$ \begin{array}{l}
\{ \rho : G \rightarrow {\rm SL}_2(R) \, | \, \mbox{deformation of}\; \overline{\rho} \; \mbox{over}\; R\}/\approx  \\
\;\;\;\;\;\;\;\;\;\;\;\;\;\;\;\;\;\;  \longrightarrow \{ T : G \rightarrow R \,  | \, \mbox{deformation of} \; {\rm tr}(\overline{\rho}) \; \mbox{over} \; R \}.
\end{array}
$$
\vspace{.1cm}

Now, by Theorem 1.1.1, there exists the universal deformation $({\boldsymbol R}_{\overline{T}}, {\boldsymbol T})$ of a pseudo-representation $\overline{T} = {\rm tr}(\overline{\rho})$.  By Theorem 1.2.1, we have a deformation $\mbox{\boldmath $\rho$} : G \rightarrow {\rm SL}_2({\boldsymbol R}_{\overline{T}})$ of $\overline{\rho}$ such that ${\rm tr}(\mbox{\boldmath $\rho$}) = {\boldsymbol T}$. Then we can verify  that $({\boldsymbol R}_{\overline{T}}, \mbox{\boldmath $\rho$})$ satisfies the desired property of the universal deformation of $\overline{\rho}$.  \\
\\
{\bf Theorem 1.2.2} ([MTTU; Theorem 2.2.2]). {\em Let $\overline{\rho} : G \rightarrow {\rm SL}_2(k)$ be an absolutely irreducible representation. Then there exists the universal deformation $({\boldsymbol R}_{\overline{\rho}}, {\boldsymbol \rho})$ of $\overline{\rho}$, where ${\boldsymbol R}_{\overline{\rho}}$ is given as ${\boldsymbol R}_{\overline{T}}$ for $\overline{T}:= {\rm tr}(\overline{\rho})$ in Theorem 1.1.1.}\\
\\
1.3. {\em Obstructions.} We recall basic facts on a presentation of a complete local $\mathcal{O}$-algebra  and the obstruction for the deformation problem. For $R \in \frak{CL}_{\cal O}$, we define the {\em relative cotangent space} $\frak{t}_{R/\mathcal{O}}^*$ of $R$ by the $k$-vector space  $\frak{m}_{R}/(\frak{m}_R^2 + \frak{m}_{\mathcal{O}}R)$ and the {\em relative tangent space} $\frak{t}_{R/\mathcal{O}}$ of $R$ by the dual $k$-vector space of $\frak{t}_{R/\mathcal{O}}^*$. We note that they are same as the cotangent and tangent spaces of $R/\frak{m}_{\mathcal{O}}R = R \otimes_{\mathcal{O}} k$, respectively. The following lemma is a well-known fact which can be proved using Nakayama's lemma (cf. [Ti; Lemma 5.1]).\\
\\
{\bf Lemma 1.3.1.} {\em Let $ d := \dim_k \frak{t}_{R/\mathcal{O}}$ and assume $d < \infty$.  Let $x_1,\dots, x_d$ be elements of $R$ whose images in $R \otimes_{\cal O} k$ form a system of parameters of  the local $k$-algebra $R \otimes_{\mathcal{O}} k$. Then  there is a surjective $\mathcal{O}$-algebra homomorphism  $$\eta : {\cal O}[[X_1, \dots , X_d]] \longrightarrow R$$
in $\frak{CL}_{\cal O}$ such that $\eta(X_i) = x_i$ for  $1\leq i \leq d$.}\\
\\
Let ${\rm Ad}(\overline{\rho})$ be the $k$-vector space ${\rm sl}_2(k) := \{ X \in {\rm M}_2(k)\, | \,{\rm tr}(X) = 0 \}$ on which $G$ acts by $g.X := \overline{\rho}(g)X\overline{\rho}(g)^{-1}$ for $g \in G$ and $X \in {\rm sl}_2(k)$. It is well-known  ([Ma2; 1.6]) that there is a canonical isomorphism between the relative cotangent space $ \frak{t}_{{\boldsymbol R}_{\overline{\rho}}/\mathcal{O}}^*$ and the 1st group cohomology $H^1(G,{\rm Ad}(\overline{\rho})).$
We say that the deformation problem for $\overline{\rho}$ is {\em unobstructed} if the 2nd cohomology $H^2(G, {\rm Ad}(\overline{\rho}))$ vanishes. The following proposition is also well-known.\\
\\
{\bf Proposition 1.3.2.} ([Ma2; 1.6, Proposition 2]). {\em Suppose that the deformation problem for $\overline{\rho}$ is unobstructed and $\dim_k H^1(G, {\rm Ad}(\overline{\rho})) < \infty$. Then the map $\eta$ in Lemma 1.3.1 with $R = {\boldsymbol R}_{\overline{\rho}}$ is isomorphic}
$$\eta : {\cal O}[[X_1, \dots , X_d]] \stackrel{\sim}{\longrightarrow} {\boldsymbol R}_{\overline{\rho}}.$$
\\
In this paper, we are interested in the case that $G$ is a knot group, namely, the fundamental group of the complement of a knot in the $3$-sphere $S^3$. We note that the deformation problem is not unobstructed in general
for a knot group representation $\overline{\rho}$, as shown in Subsection 2.3.\\

\begin{center}
{\bf 2. Character schemes}
\end{center}

In this section, we show the relation between the universal deformation ring in Section 1 and the character scheme of ${\rm SL}_2$-representations. 

In Subsection 2.1, we recall the constructions and some facts concerning the character scheme and the skein algebra over $\mathbb{Z}$, and then describe their relation.
For the details on the materials, we consult [CS], [LM, Chapter 1], [Na] and [Sa]. In Subsection 2.2, via the skein algebra, we show that the universal deformation ring may be seen as an infinitesimal deformation of the character algebra. In Subsection 2.3, we show that the deformation problem is not unobstructed for a knot group in general, using Thurston's result on the character variety. 
\\
\\
2.1. {\em Character schemes and skein algebras over $\mathbb{Z}$.} Let $G$ be a group.  Let ${\cal F}$ be the functor from the category $\frak{Com.Ring}$ of commutative rings with identity to the category of sets defined by
$$ {\cal F}(A) := \{ G \rightarrow {\rm SL}_2(A)\, | \, \mbox{representation} \}$$
for $A \in \frak{Com.Ring}$. The functor ${\cal F}$ is represented by a pair $({\cal A}(G), \sigma_G)$, where ${\cal A}(G) \in \frak{Com.Ring}$ and $\sigma_G : G \rightarrow {\rm SL}_2({\cal A}(G))$ is a representation, which satisfies the following  universal property:  ``For any $A \in \frak{Com.Ring}$ and a representation $\rho : G \rightarrow {\rm SL}_2(A)$, there is a unique morphism $\psi : {\cal A}(G) \rightarrow A$ in $\frak{Com.Ring}$ such that $\psi \circ \sigma_G = \rho $." 
 Thus the correspondence $\psi \mapsto \psi \circ \sigma_G$ gives the bijection
$$ {\rm Hom}_{\frak{Com.Ring}}({\cal A}(G), A) \simeq  \{ G \rightarrow {\rm SL}_2(A)\, | \, \mbox{representation} \}. $$
By the universal property, the pair $({\cal A}(G), \sigma_G)$ is unique (if exists) up to isomorphism.  We call ${\cal A}(G)$ the {\em universal representation algebra} over $\mathbb{Z}$ and $\sigma_G :  G \rightarrow {\rm SL}_2({\cal A}(G))$ the {\em universal representation}.  The pair $({\cal A}(G), \sigma_G)$ is constructed as follows. Let $X(g) = (X_{ij}(g))_{1 \leq i,j  \leq 2}$ be $2 \times 2$ matrix 
whose entries $X_{ij}(g)$'s  are variables indexed by $1\leq i,j\leq 2$ and $g \in G$. Then ${\cal A}(G)$ is given as
$$ {\cal A}(G) = \mathbb{Z}[ X_{ij}(g) \;  (1\leq i,j \leq 2; g \in G) ]/J,$$
where $J$ is the ideal of the polynomial ring $\mathbb{Z}[ X_{ij}(g) \; (1\leq i,j \leq 2; g \in G) ]$ generated by
$$ X_{ij}(e) - \delta_{ij}, \; X_{ij}(g_1 g_2) - \sum_{k=1}^2 X_{ik}(g_1)X_{kj}(g_2), \; \det(X(g)) - 1$$
for $1 \leq i, j \leq 2$ and $g \in G$, and  the representation $\sigma_G : G \rightarrow {\rm SL}_2({\cal A}(G))$ is given by
$$ \sigma_G(g) := X(g) \; \mbox{mod} \; J  \;\;\;\; (g \in G).$$
We note that when $G$ is presented by finitely many generators $g_1,\dots, g_n$ subject to the relations $r_l = 1 \;(l \in L)$, ${\cal A}(G)$ is given by 
$$ {\cal A}(G) = \mathbb{Z}[X_{ij}(g_h) \; (1\leq h \leq n,  1\leq i,j\leq 2)]/J' $$ 
for the ideal $J'$ generated by 
$$ r_l(X(g_1), \dots , X(g_n))_{ij} -\delta_{ij}, \;  \; {\rm det}(X(g_h)) - 1,  $$
 where $1\leq i,j\leq 2, l \in L,  1\leq h \leq n$ and $r_l(X(g_1), \dots , X(g_n))_{ij}$ denotes the $(i,j)$-entry of $r_l(X(g_1), \dots , X(g_n))$. The universal representation $\sigma_G$ is given by 
$$\sigma_G(g_h) = X(g_h) \; \mbox{mod} \;  J'   \;\;  (1\leq h \leq n).$$
So ${\cal A}(G)$ is a finitely generated algebra over $\mathbb{Z}$ if $G$ is a finitely generated group. We denote by ${\cal R}(G)$ the affine scheme ${\rm Spec}({\cal A}(G) )$  and call it the {\it  representation scheme} of $G$ over $\mathbb{Z}$. So $A$-rational points of ${\cal R}(G)$ corresponds bijectively to representations $G \rightarrow {\rm SL}_2(A)$ for any $A \in \frak{Com.Ring}$. For $\frak{p} \in {\cal R}(G)$, we let  $\rho_\frak{p} := \psi_\frak{p} \circ \sigma_G : G \rightarrow {\rm SL}_2({\cal A}(G) /\frak{p})$ be the corresponding representation,  where $\psi_\frak{p} : {\cal A}(G)  \rightarrow {\cal A}(G) /\frak{p}$ is the natural homomorphism. 

We  say that a representation $\rho : G \rightarrow {\rm SL}_2(A)$ with $A \in \frak{Com.Ring}$ is  {\em absolutely irreducible}  if the composite of $\rho$ with the natural map ${\rm SL}_2(A) \rightarrow {\rm SL}_2(k(\frak{p}))$ is absolutely irreducible over the residue field $k(\frak{p}) = A_{\frak{p}}/\frak{p}A_{\frak{p}}$ for any $\frak{p} \in {\rm Spec}(A)$. 

Let ${\rm PGL}_2$ be the group scheme over $\mathbb{Z}$ whose coordinate ring $A({\rm PGL}_2)$ is the subring of the graded ring $\mathbb{Z}[Y_{ij} \; (1\leq i,j\leq 2)]_{{\rm det}(Y)}$  consisting of homogeneous elements of degree $0$, where the degree of $Y_{ij}$ is 1. The adjoint action $ {\rm Ad} : {\cal R}(G) \times {\rm PGL}_2 \rightarrow {\cal R}(G)$ is given by  the dual action
$$ {\rm Ad}^{*} : {\cal A}(G) \longrightarrow {\cal A}(G)\otimes_{\mathbb{Z}} A({\rm PGL}_2); \;\; X_{ij}(g) \mapsto (YX(g)Y^{-1})_{ij} \otimes Y_{kl},$$
where $Y = (Y_{ij})_{1\leq i,j\leq 2}$ and $(YX(g)Y^{-1})_{ij}$ denotes the $(i,j)$-entry of $YX(g)Y^{-1}$. 
Let ${\cal B}(G)$ be the invariant subalgebra of ${\cal A}(G)$ under this action of ${\rm PGL}_2$
$$ \begin{array}{ll} {\cal B}(G) & := {\cal A}(G)^{{\rm PGL}_2}\\
                                                    & := \{ x \in {\cal A}(G) \, | \, {\rm Ad}^*(x) = x \otimes 1 \}.
  \end{array}                                        
$$
We call ${\cal B}(G)$ the {\em character algebra} of $G$ over $\mathbb{Z}$. We  denote by ${\cal X}(G)$ the affine scheme ${\rm Spec}({\cal B}(G))$ and call it  the {\em character scheme} of $G$ over $\mathbb{Z}$. The natural inclusion 
$$\iota : {\cal B}(G) \hookrightarrow {\cal A}(G)$$
 induces a morphism of schemes 
$$\iota^{\#} : {\cal R}(G) \longrightarrow {\cal X}(G).$$
 We denote the image of $\frak{p} (= \rho_{\frak{p}}) \in {\cal R}(G)$ in ${\cal X}(G)$ under  $\iota^{\#}$  by $[\frak{p}] (=[\rho_{\frak{p}}])$. 

 According to [PS; Definition 2.5] and [Sa; 3.1],  we define  the {\em skein algebra} ${\cal C}(G)$ over $\mathbb{Z}$ by
 $$ {\cal C}(G) := \mathbb{Z}[t_g \; (g \in G)]/I,$$
where $I$ is the ideal of the polynomial ring $\mathbb{Z}[t_g  \; (g \in \Pi)]$ generated by the polynomials of the form
$$t_e -2,\;\;  t_{g_1}t_{g_2} - t_{g_1g_2} - t_{g_1^{-1}g_2} \;\;\;\; (g_1, g_2 \in G). $$
We note that ${\cal C}(G)$ is a finitely generated algebra over $\mathbb{Z}$ if  $G$ is a finitely generated group ([Sa; 3.2]). We denote by ${\cal S}(G)$ the affine scheme ${\rm Spec}({\cal C}(G))$  and call it the {\em  skein scheme} of $G$  over $\mathbb{Z}$. 

Since ${\rm tr}(\sigma_G(g)) \; (g \in G)$ is invariant under the adjoint action of ${\rm PGL}_2$ and we have the formula
$$ {\rm tr}(\sigma_G(g_1)) {\rm tr}(\sigma_G(g_2)) - {\rm tr}(\sigma_G(g_1g_2)) - {\rm tr}(\sigma_G(g_1^{-1}g_2)) = 0$$
for $g_1, g_2 \in G$, which is derived from the Cayley-Hamilton relation, we obtain a $\mathbb{Z}$-algebra homomorphism 
$$ \tau : {\cal C}(G) \longrightarrow {\cal B}(G)$$
defined by 
$$\tau(t_g) := {\rm tr}(\sigma_G(g)) \;\;\;\;\; (g \in G).$$
It induces the morphism of schemes
$$ \tau^{\#} : {\cal X}(G) \longrightarrow {\cal S}(G).$$
We set 
$$\varphi := \iota \circ \tau : {\cal C}(G) \longrightarrow {\cal A}(G) $$
 so that we have the morphism of schemes
$$ \varphi^{\#}  = \tau^{\#} \circ \iota^{\#} : {\cal R}(G)  \longrightarrow  {\cal S}(G).$$

Now we define the {\it discriminant ideal} $\Delta(G)$ of ${\cal C}(G)$ by the ideal generated by the images of the elements in $\mathbb{Z}[t_g \; (g \in \pi)]$ of the form
$$ \Delta(g_1,g_2) := t_{g_1g_2g_1^{-1}g_2^{-1}} -2 = t_{g_1}^2 + t_{g_2}^2 + t_{g_1g_2}^2 - t_{g_1}t_{g_2}t_{g_1g_2} - 4 \;\; (g_1,g_2 \in G),$$
and the {\it discriminant subscheme} by $V(\Delta(G)) = {\rm Spec}({\cal C}(G)/\Delta(G))$.  We define the open subschemes ${\cal S}(G)_{\rm a.i}$, ${\cal X}(G)_{\rm a.i}$ and ${\cal R}(G)_{\rm a.i}$ of ${\cal S}(G)$, ${\cal X}(G)$ and ${\cal R}(G)$, respectively,  by
$$\begin{array}{l}  {\cal S}(G)_{\rm a.i} :=  {\cal S}(G) \setminus V(\Delta(G)), \\
 {\cal X}(G)_{\rm a.i} := {\cal X}(G) \setminus (\tau^{\#})^{-1}(V(\Delta(G))),\\
 {\cal R}(G)_{\rm a.i} := {\cal R}(G) \setminus (\varphi^{\#})^{-1}(V(\Delta(G))).
\end{array}$$
The following theorem, due to Kyoji Saito, is fundamental for our purpose.
\\
\\
{\bf Theorem 2.1.1} ([Sa; 4.2, 4.3], [Na; Corollary 6.8]). (1) {\em For $\frak{p} \in {\cal R}(G)$, $\rho_{\frak{p}}$ is absolutely irreducible if and only if $\frak{p} \in  {\cal R}(G)_{\rm a.i}.$}\\
(2) {\em The restriction of $\varphi^{\#}$ to ${\cal R}(G)_{\rm a.i}$ 
$$ \varphi^{\#}_{\rm a.i} : {\cal R}(G)_{\rm a.i} \longrightarrow {\cal S}(G)_{\rm a.i}$$
is a principal ${\rm PGL}_2$-bundle.}\\
(3) {\em The restriction of $\tau^{\#}$ to ${\cal X}(G)_{\rm a.i}$ is an isomorphism}
$$ \tau^{\#}_{\rm a.i} : {\cal X}(G)_{\rm a.i} \stackrel{\sim}{\longrightarrow} {\cal S}(G)_{\rm a.i}.$$
\\
By virtue of Theorem 2.1.1 (1), we call ${\cal S}(G)_{\rm a.i}$, ${\cal X}(G)_{\rm a.i}$ and ${\cal R}(G)_{\rm a.i}$ the {\it absolutely irreducible part} of ${\cal S}(G)$, $ {\cal X}(G)$ and ${\cal R}(G)$, respectively. We note that ${\cal X}(G)_{\rm a.i} (\simeq {\cal S}(G)_{\rm a.i})$ represents the functor $\overline{\cal F}$ from the category $\frak{Sch}$ of schemes to the category of sets, which associates to a scheme $X$ the set of isomorphism classes of absolutely irreducible representations $G \rightarrow {\rm SL}_2(\Gamma(X, {\cal O}_X))$:
$$ \overline{\cal F}(X) := \{ G \rightarrow {\rm SL}_2(\Gamma(X, {\cal O}_X))\, | \, \mbox{absolutely irreducible representation} \}/\sim.$$
 Since $ \varphi^{\#}_{\rm a.i}$ and $\tau^{\#}_{\rm a.i}$ are defined over $\mathbb{Z}$,   they induces maps on $A$-rational points for $A \in \frak{Com.Ring}$:
$$\begin{array}{l} \varphi^{\#}_{\rm a.i}(A) : {\cal R}(G)_{\rm a.i}(A) \longrightarrow {\cal S}(G)_{\rm a.i}(A),\\
\tau^{\#}_{\rm a.i}(A) : {\cal X}(G)_{\rm a.i}(A) \stackrel{\sim}{\longrightarrow} {\cal S}(G)_{\rm a.i}(A).
\end{array}
$$
By Theorem 2.1.1 (3), we have the following\\
\\
{\bf Corollary 2.1.2.} {\it Let  $\overline{\rho} : G \rightarrow {\rm SL}_2(k)$ be an absolutely irreducible representation  over a field $k$ so that $\overline{\rho} \in {\cal R}(G)_{\rm a.i}(k)$. Let $[\overline{\rho}] \in {\cal X}(G)_{\rm a.i}(k)$ also denote the corresponding prime ideal of ${\cal B}(G)$. Then the morphism $\tau$ induces an isomorphism of local rings}:
$$ {\cal C}(G)_{\tau^{\#}([\overline{\rho}])} \simeq {\cal B}(G)_{[\overline{\rho}]}.$$
\\
The following proposition can be proved by using the vanishing of the Galois cohomology $H^1(k,{\rm PGL}_2(\overline{k})) = 1$ for a field $k$ whose Brauer group ${\rm Br}(k) = 0$ ([Se3; III, 2.2]) and Skolem-Noether theorem. For example, when $k$ is a finite field or an algebraically closed field, ${\rm Br}(k) = 0$.\\
\\
{\bf Proposition 2.1.3} ([Fu; Lemma 3.3.1], [Ha; Proposition 2.2.27]). {\em Let $k$ be a field whose Brauer group ${\rm Br}(k) = 0$. Then $\varphi^{\#}_{\rm a.i}$ induces the following bijection on $k$-rational points}:
$$ \varphi^{\#}_{\rm a.i}(k) : \; {\cal R}(G)_{\rm a.i}(k)/{\rm PGL}_2(k) \stackrel{\sim}{\longrightarrow} {\cal S}(G)_{\rm a.i}(k).$$
\\
2.2. {\em The relation between the universal deformation ring and the character scheme.}  Let $k$ be  a perfect  field with ${\rm char}(k) \neq 2$ and let ${\cal O}$ be a discrete valuation ring  with residue field $k$. Let $\overline{\rho} : G \rightarrow {\rm SL}_2(k)$ be an absolutely  irreducible representation and let $\overline{T} : G \rightarrow k$  be a pseudo-${\rm SL}_2$-representation over $k$ given by the character ${\rm tr}(\overline{\rho})$. Let ${\boldsymbol R}_{\overline{\rho}} (= {\boldsymbol R}_{\overline{T}})$ be the universal deformation ring of $\overline{\rho}$ (or $\overline{T}$) as in Theorem 1.2.2. Recall that ${\boldsymbol R}_{\overline{T}}$ is a complete local ${\cal O}$-algebra whose  residue field  is  $k$.  On the other hand, let ${\cal B}(G)$ and ${\cal S}(G)$ be the character algebra and skein algebra of $G$ over $\mathbb{Z}$, respectively.  We set
$$\begin{array}{ll} {\cal B}(G)_k := {\cal B}(G) \otimes_{\mathbb{Z}} k, & {\cal X}(G)_k := {\rm Spec}({\cal B}(G)_k) = {\cal X}(G)\otimes_{\mathbb Z} k,\\
{\cal C}(G)_k := {\cal C}(G) \otimes_{\mathbb{Z}} k, & {\cal S}(G)_k := {\rm Spec}({\cal C}(G)_k) = {\cal S}(G)\otimes_{\mathbb Z} k. 
\end{array}$$
We also denote by ${\cal X}(G)_k^{\rm a.i}$ and ${\cal S}(G)_k^{\rm a.i}$ the absolutely irreducible part of ${\cal X}(G)_k$ and ${\cal S}(G)_k$, respectively. By Theorem 2.1.1 (3), we have ${\cal X}(G)_k^{\rm a.i} \simeq {\cal S}(G)_k^{\rm a.i}$. The following theorem tells us that the universal deformation ring ${\boldsymbol R}_{\overline{\rho}}$ may be seen as an infinitesimal deformation of the character $k$-algebra ${\cal B}(G)_k$ at $[\overline{\rho}]$.\\
\\
{\bf Theorem 2.2.1.} {\em Let $[\overline{\rho}]$ denote the maximal ideal of ${\cal B}(G)_k$ corresponding to the representation $\overline{\rho}$. We then have an isomorphism of $k$-algebras}
$$ {\boldsymbol R}_{\overline{\rho}}\otimes_{\cal O} k  \; \simeq \; ({\cal B}(G)_k)_{[\overline{\rho}]}^{\wedge},$$
{\em where $({\cal B}(G)_k)_{[\overline{\rho}]}^{\wedge}$ denotes the $[\overline{\rho}]$-adic completion of ${\cal B}(G)_k$. }\\
\\
{\em Proof.} By the construction of  ${\boldsymbol R}_{\overline{T}}$ in Theorem 1.1.1, we have
$$ {\boldsymbol R}_{\overline{\rho}} = {\cal O}[[ X_g \; (g \in G)]]/{\cal I},$$
where ${\cal I}$ is the ideal of the power series ring ${\cal O}[[ X_g \; (g \in G)]]$ generated by elements of the form: setting $T_g := X_g + \lambda(\overline{T}(g))$ ($\lambda$ : the Teichm\"{u}ller lift),
\vspace{.1cm}\\
(1) $T_e -2$,\\
(2) $T_{g_1g_2} - T_{g_2g_1}$,\\
(3) $T_{g_1}T_{g_2}T_{g_3} + T_{g_1g_2g_3} + T_{g_1g_3g_2} -T_{g_1g_2}T_{g_3} - T_{g_2g_3}T_{g_1} - T_{g_1g_3}T_{g_2}$,\\
(4) $T_g^2 - T_{g^2} -2,$
\vspace{.1cm}\\
where $g, g_1, g_2, g_3 \in G$. 

On the other hand, let $\psi : {\cal B}(G)_k \rightarrow k$ be the morphism in $\frak{Com.Ring}$ corresponding to $[\overline{\rho}] \in {\cal X}(G)_{\rm a.i}(k)$. Since $ \psi({\rm tr}(\sigma_G(g))) = {\rm tr}(\overline{\rho}(g)) = \overline{T}(g)$ for $g \in G$, the maximal ideal $[\overline{\rho}] = {\rm Ker}(\psi)$ of ${\cal B}(G)_k$ corresponds to the maximal ideal $(t_g - \overline{T}(g)\; (g \in G))$ of ${\cal C}(G)_k$.  Therefore  Corollary 2.1.2 yields
$$ ({\cal B}(G)_k)_{[\overline{\rho}]}^{\wedge} \simeq  k[[ x_g \; (g \in G)]]/I^{\wedge},$$
where $x_g := t_g - \overline{T}(g) \; (g \in G)$ and $I^{\wedge}$ is the ideal of the power series ring $k[[ x_g \; (g \in G)]]$ generated by  elements of the form 
$$t_e -2,\;  t_{g_1}t_{g_2} - t_{g_1g_2} - t_{g_1^{-1}g_2}\;  (g_1, g_2 \in  G).$$
So, in order to show that the correspondence  $x_g  \mapsto X_g  \otimes 1$ gives the desired isomorphism $({\cal B}(G)_k)_{[\overline{\rho}]}^{\wedge} \simeq  {\boldsymbol R}_{\overline{\rho}} \otimes_{\cal O} k$, it suffices to show the following\\
\\
{\bf Lemma 2.2.2.} {\em Let $T$ be a function on $G$ with values in an integral domain whose characteristic is not $2$. Let ${\rm (P)}$ be the relations given by}
\vspace{.12cm}\\
(P1) $T(1) = 2,$\\
(P2) $T(g_1g_2) = T(g_2g_1)$,\\
(P3) $T(g_1)T(g_2)T(g_3)+T(g_1g_2g_3)+T(g_1g_3g_2)-T(g_1g_2)T(g_3)-T(g_2g_3)T(g_1)-T(g_1g_3)T(g_2) =0 $,\\
(P4) $T(g)^2 - T(g^2) = 2$,
\vspace{.12cm}\\
{\em and let ${\rm (C)}$ be the relations given by}
\vspace{.12cm}\\
(C1) $T(1) = 2$,\\
(C2) $T(g_1)T(g_2) = T(g_1g_2) + T(g_1^{-1}g_2)$,
\vspace{.12cm}\\
{\em where $g, g_1, g_2, g_3$ are any elements in $G$. 

Then ${\rm (P)}$ and ${\rm (C)}$ are equivalent.}\\
\\
{\em Proof of} Lemma 2.2.2.  (P) $\Rightarrow$ (C): Letting $g_2 = g_1$ in (P3), we have
$$T(g_1)^2T(g_3)-T(g_1^2)T(g_3)+T(g_1^2g_3)+T(g_1g_3g_1)-2T(g_1g_3)T(g_1)=0.$$
Using (P2) and (P4), we have
$$2(T(g_3)+T(g_1^2g_3)-T(g_1g_3)T(g_1))=0.$$
Letting $g_3$ be replaced by $g_1^{-1}g_2$ in the above equation and noting $T$ has the value in an integral domain whose characteristic is not $2$, we obtain (C2).

(C) $\Rightarrow$ (P). Letting $g_2=1$ in (C2) and using (C1), we have 
$$ T(g) = T(g^{-1}) \; \mbox{for any}\; g \in G. $$
Exchanging $g_1$ and $g_2$ in (C2) each other and using the above relation, we have
$$ T(g_2)T(g_1) = T(g_2g_1) + T(g_2^{-1}g_1) = T(g_2g_1) + T(g_1^{-1}g_2)$$
and hence we obtain (P2). Next letting $g_1$ be replaced by $g_1g_3$ in (C2), we have
$$ -T(g_1g_3)T(g_2) + T(g_1g_3g_2) + T(g_3^{-1}g_1^{-1}g_2) = 0, \leqno{(2.2.2.1)} $$
and letting $g_2$ be replaced by $g_2g_3$ in (C2), we have
$$ -T(g_1)T(g_2g_3) + T(g_1g_2g_3) + T(g_1^{-1}g_2g_3) = 0. \leqno{(2.2.2.2)}$$
By (C2), we have
$$ \begin{array}{ll}
T(g_3^{-1}g_1^{-1}g_2) & = T(g_3)T(g_1^{-1}g_2) - T(g_3g_1^{-1}g_2) \\
 & = T(g_3)T(g_1)T(g_2)-T(g_1g_2)T(g_3) - T(g_3g_1^{-1}g_2).
\end{array}$$
Hence, using (P2) proved already, we have
$$ \begin{array}{ll}
T(g_3^{-1}g_1^{-1}g_2) + T(g_1^{-1}g_2g_3) & = T(g_1)T(g_2)T(g_3)-T(g_1g_2)T(g_3) \\
                                                                           & \;\;\;\;\;\;\;\; -T(g_3g_1^{-1}g_2)+T(g_1^{-1}g_2g_3) \\
& = T(g_1)T(g_2)T(g_3)-T(g_1g_2)T(g_3). 
\end{array} \leqno{(2.2.2.3)} $$
Summing up (2.2.2.1) and (2.2.2.2) together with (2.2.2.3), we obtain (P3). Finally putting $g_1 = g_2$ in (C2) and using (C1), we obtain (P4). $\;\;\Box$\\

By Lemma 1.3.1 and Theorem 2.2.1, we have the following\\
\\
{\bf Corollary 2.2.3.} {\em Let $[\overline{\rho}]$ be a regular $k$-rational point of ${\cal X}(G)_k^{\rm a.i}$. Let $d$ be the dimension of the irreducible component of ${\cal X}(G)_k^{\rm a.i}$ containing $[\overline{\rho}]$ so that 
$({\cal B}(G)_k)_{[\overline{\rho}]}^{\wedge}$ is a power series ring over $k$ on a regular system of parameters  $z_1,\dots, z_d$. Let $x_1, \dots , x_d$ be elements of ${\boldsymbol R}_{\overline{\rho}}$ such that the image of $x_i$ in ${\boldsymbol R}_{\overline{\rho}} \otimes_{\cal O} k \simeq  ({\cal B}(G)_k)_{[\overline{\rho}]}^{\wedge}$ is $z_i$ for $1\leq i \leq d$. Then  there is  a surjective ${\cal O}$-algebra homomorphism 
$$ \eta : {\cal O}[[X_1,\dots, X_d]] \longrightarrow {\boldsymbol R}_{\overline{T}}$$
 in $\frak{CL}_{\cal O}$ such that $\eta(X_i) = x_i$ for $1\leq i \leq d$.}\\
 \\
 By Corollary 2.2.3, we obtain the following criterion which determines the universal deformations for many examples. See Section 4. \\
 \\
 {\bf Theorem 2.2.4.} {\em Let notations and assumptions be as in Corollary 2.2.3. We suppose that there are $g_1, \dots , g_d \in G$ such that $z_i = t_i - {\rm tr}(\overline{\rho}(g_i))$ for $1\leq i \leq d$, where $t_i$ denotes a variable corresponding to the regular function ${\rm tr}(\sigma_G(g_i))$.  Choose $\alpha_i  \in {\cal O}$ such that $\alpha_i \; \mbox{mod}\; \frak{m}_{\cal O} = {\rm tr}(\overline{\rho}(g_i))$ for $1\leq i \leq d$ and suppose that  $\rho : G \rightarrow {\rm SL}_2({\cal O}[[t_1 - \alpha_1,\dots , t_d - \alpha_d]])$ is a deformation of $\overline{\rho}$ satisfying
 $$ {\rm tr}(\rho(g_i)) = t_i \;\;\;\;\; (1\leq i \leq d).$$
Then $({\cal O}[[t_1 - \alpha_1, \dots , t_d - \alpha_d]], \rho)$ is the universal deformation of $\overline{\rho}$.}\\
\\ 
{\em Proof.} By the universal property of $({\boldsymbol R}_{\overline{\rho}}, {\boldsymbol \rho})$, there is a morphism 
$$ \psi : {\boldsymbol R}_{\overline{\rho}} \longrightarrow {\cal O}[[t_1 - \alpha_1, \dots , t_d-\alpha_d]]$$
in $\frak{CL}_{\cal O}$ such that $\psi \circ {\boldsymbol \rho} \approx \rho$. Hence we have
$$ \psi({\rm tr}({\boldsymbol \rho}(g_i))) = {\rm tr}(\rho(g_i)) = t_i, \;\; 1\leq i \leq d. \leqno{(2.2.4.1)} $$
By  Corollary 2.2.3,  there is a  surjective morphism 
$$\eta : {\cal O}[[X_1, \dots , X_d ]] \longrightarrow {\boldsymbol R}_{\overline{\rho}}$$
 in $\frak{CL}_{\cal O}$ such that $\eta(X_i) = {\rm tr}({\boldsymbol \rho}(g_i)) - \alpha_i$ for $1\leq i \leq d$. Since $\psi \circ \eta : {\cal O}[[X_1, \dots , X_d]] \rightarrow {\cal O}[[t_1-\alpha_1, \dots , t_d - \alpha_d]]$
is a morphism in $\frak{CL}_{\cal O}$ and satisfies, by (2.2.4.1), 
$$\psi \circ \eta (X_i) = t_i - \alpha_i \;\;\;\;\; (1\leq i \leq d),$$ 
$\psi \circ \eta$ is an isomorphism in $\frak{CL}_{\cal O}$. Since $\eta$ is surjective,
$\eta$ must be isomorphic and so is $\psi$. $\;\; \Box$\\
\\
2.3. {\em The case of a knot group}.  Let $K$ be a knot in the $3$-sphere $S^3$ and  let $E_K$  denote the knot complement $S^3 \setminus K$. Let $G_K$ denote the knot group of $K$,  $G_K := \pi_1(E_K)$. It is well known that  
$G_K$ has the following presentation of deficiency one (for example, the Wirtinger presentation):
$$
G_K=\left< g_1, \dots, g_n \; | \; r_1 = \cdots = r_{n-1} = 1 \right>. \leqno{(2.3.1)}
$$
Let $k$ be a field with ${\rm char}(k) \neq 2$. Let  $\overline{\rho} : G_K \rightarrow {\rm SL}_2(k)$ be an absolutely irreducible representation and let ${\boldsymbol R}_{\overline{\rho}}$ be the universal deformation ring as in Theorem 1.2.2. 
Since the character variety ${\cal X}(G_K)_k$ of a knot group $G_K$ over a field $k$ has been extensively studied (see [CS], [Le], [Ha] etc),  we can determine ${\boldsymbol R}_{\overline{\rho}} \otimes_{\cal O} k$ by Theorem 2.2.1 and even ${\boldsymbol  R}_{\overline{\rho}}$ by Theorem 2.2.4 for some knots $K$. In fact, in [MTTU], we determined ${\boldsymbol R}_{\overline{\rho}}$ for a certain Riley-type representations $\overline{\rho}$ for a 2-bridge knots  $K$. See also Section 4 for other examples. 

It is a delicate problem, however, to determine the universal deformation ring ${\boldsymbol R}_{\overline{\rho}}$ for a  knot group representation $\overline{\rho}$ in general, since the deformation problem for $\overline{\rho}$ is not unobstructed in general for a knot group $G_K$, as the following theorem shows.\\
\\
{\bf Theorem 2.3.2.} {\em We suppose that  $\rho : G_K \rightarrow {\rm SL}_2(\mathbb{C})$ is an irreducible representation and that  there is a subring $A$ of a finite algebraic number field $F$ and a finite prime $\frak{p}$ of $F$ such that $A$ is $\frak{p}$-integral and the image of $\rho$ is contained in ${\rm SL}_2(A)$. Set $k := A/\frak{p}$ and $\overline{\rho} := \rho \; \mbox{mod}\; \frak{p} : G_K \rightarrow {\rm SL}_2(k)$. Then we have}
$$ H^2(G_K, {\rm Ad}(\overline{\rho})) \neq 0.$$
\\
We note that the assumption in Theorem 2.3.2 is satisfied, for instance, when $K$ is a hyperbolic knot and $\rho$ is the holonomy representation attached to a hyperbolic structure on $E_K$ such that the completion is a closed or a cone 3-manifold.
For the proof of Theorem 2.3.2, we recall the following lemma, which is a special case of a more general result, due to Thurston, for 3-manifolds.\\
\\
{\bf Lemma 2.3.3} ([CS; Proposition 3.2.1]). {\em For an irreducible representation $\rho : G_K \rightarrow {\rm SL}_2(\mathbb{C})$, the irreducible component of ${\cal X}(G_K)_{\mathbb{C}}$ containing $[\rho]$ has the dimension greater than $0$.}\\
\\
{\em Proof of Theorem 2.3.2.} Let $W$ be the CW complex attached to the presentation (2.3.1). We recall herewith the construction of $W$:\\
$\bullet$ We prepare  $0$-cell $b^*$, $1$-cells $g_1^*, \dots , g_n^*$, where each $g_i^*$  corresponds to the generator $g_i$, $2$-cells $r_1^*, \dots , r_{n-1}^*$, where each $r_j^*$ corresponds to the relator $r_j$.\\
$\bullet$ We attach each $1$-cell $g_i^*$ to the $0$-cell $b^*$ so that we obtain a bouquet.\\
$\bullet$ We attach the boundary of each $2$-cell $r_j^*$ to $1$-cells of the bouquet, according to words in $r_j$.\\
We note that the knot complement $E_K$ and the CW complex $W$ are homotopically equivalent by Whitehead's theorem, because they are  both the Eilenberg-MacLane space $K(G_K,1)$. \\

 Let  ${\rm Ad}(\rho)$ be the $A$-module ${\rm sl}_2(A)$ on which $G_K$ acts by $g.X :=  \rho(g)X\rho(g)^{-1}$ for $g \in G_K$ and $X \in {\rm sl}_2(A)$.  We let ${\rm Ad}(\rho)_{\mathbb{C}} := {\rm Ad}(\rho) \otimes_A \mathbb{C} = {\rm sl}_2(\mathbb{C})$ on which $G_K$ acts as $g \otimes {\rm id}_{\mathbb{C}}$ for $g \in G_K$.  Since the Euler characteristic of $W$ is zero, we have
$$ \begin{array}{ll}
\displaystyle{\sum_{i=0}^2 (-1)^i \dim_{\mathbb{C}} H^i(G_K,{\rm Ad}(\rho)_{\mathbb{C}})} & = \displaystyle{ \sum_{i=0}^2 (-1)^i \dim_{\mathbb{C}} C^i(W; {\rm Ad}(\rho)_{\mathbb{C}}) }\\
                                                                                                         & = \displaystyle{ 3 \sum_{i=0}^2 (-1)^i \dim_{\mathbb{C}} C^i(W;\mathbb{C}) }\\
                                                                                                         & = 0.
\end{array}
\leqno{(2.3.2.1)}
$$
Since $\rho$ is irreducible, we have $H^0(G_K,{\rm Ad}(\rho)_{\mathbb{C}}) = 0$ by Schur's lemma. So, by (2.3.2.1), we  have   
$$\dim_{\mathbb{C}} H^2(G_K,{\rm Ad}(\rho)_{\mathbb{C}})   =  \dim_{\mathbb{C}} H^1(G_K,{\rm Ad}(\rho)_{\mathbb{C}}). \leqno{(2.3.2.2)} $$  Since $H^1(G_K,{\rm Ad}(\rho)_{\mathbb{C}})$ contains the tangent space of the character variety ${\cal X}(G_K)_{\mathbb{C}}$ at $[\rho]$ ([Po; Proposition 3.5]), Lemma 2.3.3 implies  $H^1(G_K, {\rm Ad}(\rho))_{\mathbb{C}}) \neq 0$. So, by (2.3.2.2), we have   $H^2(G_K,{\rm Ad}(\rho)_{\mathbb{C}}) \neq 0$. Since $H^2(G_K, {\rm Ad}(\rho)_{\mathbb{C}}) = H^2(G_K, {\rm Ad}(\rho)_A) \otimes_A \mathbb{C}$, we have 
$$H^2(G_K, {\rm Ad}(\rho)) \neq 0.  \leqno{(2.3.2.3)}$$                                                                         

Let ${\rm Ad}(\overline{\rho}) := {\rm Ad}(\rho) \otimes_A k = {\rm sl}_2(k)$ on which $G_K$ acts as $g \otimes {\rm id}_k$ for $g \in G_K$. 
Let us  consider the differentials of cochains
$$ \begin{array}{l}
d : C^1(W; {\rm Ad}(\rho))  \longrightarrow C^2(W; {\rm Ad}(\rho)),\\
\overline{d} := d \otimes (\mbox{mod} \; \frak{p}) :  C^1(W; {\rm Ad}(\overline{\rho}))  \longrightarrow C^2(W; {\rm Ad}(\overline{\rho})).
\end{array}
$$
By (2.3.2.3), all $3n$-minors of $d$ are  zero. Therefore  all $3n$-minors of $\overline{d}$ are zero and  hence $H^2(G_K,{\rm Ad}(\overline{\rho})) \neq 0$. $\;\; \Box$.\\

\begin{center}
{\bf 3. $L$-functions associated to the  universal deformations}
\end{center}

In this section, we study the twisted knot module  $H_1({\boldsymbol \rho})  = H_1(E_K; {\boldsymbol \rho})$  with coefficients in the universal deformation ${\boldsymbol \rho}$ of an ${\rm SL}_2$- representation of a knot group $G_K$, and introduce the associated $L$-function $L_K({\boldsymbol \rho})$. We then formulate two problems proposed by Mazur ([Ma3]):  the torsion property of $H_1({\boldsymbol \rho})$ over the universal deformation ring ${\boldsymbol R}_{\overline{\rho}}$ (Problem 3.2.3) and the generic simplicity of the zeroes of $L_K({\boldsymbol \rho})$ (Problem 3.2.10). Our main theorem in this section (Theorem 3.2.4) gives a criterion for $H_1({\boldsymbol \rho})$  to be  finitely generated and torsion over ${\boldsymbol R}_{\overline{\rho}}$ using  a twisted Alexander invariant of $K$.  
\\
\\
3.1. {\em Fitting ideals and twisted Alexander invariants.}  Let $A$ be a Noetherian integrally closed domain. Let $M, M'$ be finitely generated $A$-modules. We say that a homomorphism $\varphi : M \rightarrow M'$ is a {\em pseudo-isomorphism} if the annihilators of ${\rm Ker}(f)$ and ${\rm Coker}(f)$ are not contained in height $1$ prime ideals of $A$. \\
\\
{\bf Lemma 3.1.1} (cf. [Se1; Lemma 5]). {\em  For any finitely generated torsion $A$-module $M$, there are positive integers $e_1, \dots , e_s$, height $1$ prime ideals $\frak{p}_1, \dots , \frak{p}_s$ of $A$ for some $s \geq 1$,  and
a pseudo-isomorphism}
$$ \varphi : M \longrightarrow \bigoplus_{i=1}^s A/\frak{p}_i^{e_i}.$$
{\em Here the set $\{ (\frak{p}_i,  e_i) \}$ is uniquely determined by $M$. If $A$ is a Noetherian factorial domain further, each prime ideal $\frak{p}_i$ of height $1$  is a principal ideal $\frak{p}_i = (f_i)$ for a prime element $f_i$ of $A$.}
\\
\\
We note that a regular local ring is a Noetherian factorial local domain   (Auslander-Buchsbaum). For example, the Iwasawa algebra ${\cal O}[[X]]$ is a 2-dimensional regular local ring, where ${\cal O}$ is a complete discrete valuation ring with ${\rm char }({\cal O}) = 0$ and finite residue field. Then it is known in Iwasawa theory ([I]) that  a height $1$ prime ideal of ${\cal O}[[X]]$  is $(\varpi)$ for a prime element $\varpi$ of ${\cal O}$ or $(f)$ for an irreducible distinguished polynomial $f \in {\cal O}[X]$, and a pseudo-isomorphism means a homomorphism with finite kernel and cokernel ([Ws; $\S$13.2]).

Let $A$ be a Noetherian factorial domain and let $M$ be a finitely generated  $A$-module. Let us take a finite presentation of $M$ over $A$
$$ A^m \stackrel{\partial}{\longrightarrow} A^n \longrightarrow M \longrightarrow 0,$$
where $\partial$ is an $n \times m$ matrix over $A$. For a non-negative integer $d$, we define  the $d$-th {\em Fitting ideal} ({\em elementary ideal}) $E_d(M)$ of $M$ to be the ideal generated by $(n-d)$ minors of $\partial$. If $d \geq n$, we let $E_d(M) := A$, and if $n-d > m$, we let  $E_d(M) := 0$. These ideals depend only on $M$ and independent of the choice of a presentation. The initial Fitting ideal $E_0(M)$ is called the {\em order ideal} of $M$. Let $\Delta_d(M)$ be the greatest common divisor of generators of $E_d(M)$, which is well defined up to multiplication by a unit of $A$. The {\em rank} of $M$ over $A$ is defined by the dimension of $M \otimes_A Q(A)$ over $Q(A)$. The following facts are well known ([Kw; 7.2], [Hi; Ch.3]).\\
\\
{\bf Lemma 3.1.2.} {\em Let $0 \rightarrow M_1 \rightarrow M_2 \rightarrow M_3 \rightarrow 0$ be an exact sequence of finitely generated $A$-modules. Then we have the followings.}\\
(1) {\em $\Delta_0(M_2) \, \dot{=} \, \Delta_0(M_1) \Delta_0(M_3).$ }\\
(2) {\em If the $A$-torsion subgroup of $M_3$ is zero and $r$ is the rank of $M_3$ over $A$, then $\Delta_d(M_2) \, \dot{=} \, \Delta_{d-r}(M_1)$.}\\
\\
For example, suppose $A$ is a principal ideal domain and $M$ is a finitely generated torsion $A$-module. Then we have $M \simeq \bigoplus_{i=1}^s A/(a_i)$ with $(a_1) \supset \cdots \supset (a_s)$, and  $E_d(M) = (a_1\cdots a_{s-d})$, $\Delta_d(M) \, \dot{= } \, a_1\cdots a_{s-d}$ for $d < s$. As another example, let $A$ be the Iwasawa algebra ${\cal O}[[X]]$ and $M$ a finitely generated torsion $A$-module. Then there is a pseudo-isomorphism $\varphi : M \rightarrow  \bigoplus_{i=1}^s  A/(f_i^{e_i})$, where $f_i$ is a prime element of ${\cal O}$ or an irreducible distinguished polynomial in ${\cal O}[X]$. If $\varphi$ is injective, in particular, if $M$ has no non-trivial finite $A$-submodule, we have $E_0(M) = (f)$, $\Delta_0(M) \, \dot{=} \,  f$, where  $f $ is the {\em Iwasawa polynomial} $\prod_{i=1}^s f_i^{e_i}$  ([MW; Appendix]). For higher Fitting ideals $E_d(M)$ for $d > 0$ in Iwasawa theory, we refer to [Ku].

Next, let $C$ be a finite connected CW complex. Let $G := \pi_1(C)$ be the fundamental group of $C$ which is supposed to have the finite presentation
$$
G=\left< g_1, \dots, g_n \; | \; r_1 = \cdots = r_{m} = 1 \right>, 
$$
where relators $r_1, \dots, r_{m}$ are words of the letters $g_1, \dots , g_n$. We suppose that there is a surjective homomorphism
$$ \alpha : G \longrightarrow \langle t \rangle \simeq \mathbb{Z}.$$
Let $A$ be a Noetherian factorial domain. We denote by the same $\alpha$ the group $A$-algebra homomorphism $A[G] \rightarrow A[ t^{\pm 1}]$, which is induced by $\alpha$. Let
$$ \rho : G \longrightarrow {\rm GL}_N(A)$$
be a representation of $G$ of degree $N$ over $A$ and let us denote by the same $\rho$ the $A$-algebra homomorphism $A[G] \rightarrow {\rm M}_N(A)$ induced by $\rho$. Then we have the tensor product representation
$$ \rho \otimes \alpha : A[G] \longrightarrow {\rm M}_N(A[t^{\pm 1}]).$$
The twisted Alexander invariant $\Delta(C,\rho; t) \in A[t^{\pm 1}]$ is defined as follows. Let $F$ be the free group on $g_1, \dots , g_n$ and let $\pi : F \rightarrow G$ be the natural homomorphism. We denote by the same $\pi$ the $A$-algebra homomorphism $A[F] \rightarrow A[G]$ induced by $\pi$. Then we have the $A$-algebra homomorphism
$$ \Phi := (\rho \otimes \alpha) \circ \pi : A[F] \longrightarrow {\rm M}_N(A[t^{\pm 1}]).$$ 
Let $\frac{\partial}{\partial g_i} : A[F] \rightarrow A[F]$ be the Fox derivative over $A$, extended from $\mathbb{Z}$ ([Fo]). Let us consider the (big) $n \times m$ matrix $P$, called the {\em twisted Alexander matrix}, whose $(i,j)$ component is defined by the $N \times N$ matrix
$$ \Phi \left( \frac{\partial r_j}{\partial g_i} \right).$$
For $1\leq i \leq n$, let $P_i$ denote the matrix obtained by deleting the $i$-th row from $P$ and we regard $P_i$ as an $(n-1)N \times mN$ matrix over $A[t^{\pm 1}]$. We note that $A[t^{\pm 1}]$ is also a Noetherian factorial domain. Let $D_i$ be the greatest common divisor of all $(n-1)N$-minors of $P_i$.
Then it is known that there is $i$ $(1\leq i \leq n)$  such that $\det(\Phi(g_i-1)) \neq 0$ and that  the ratio 
$$\Delta(C,\rho; t) := \frac{D_i}{\det \Phi(g_i-1)}  \; (\in Q(A)(t)) \leqno{(3.1.3)}$$
 is independent of such $i$'s  and is called the {\em twisted Alexander invariant} of $C$ associated to $\rho$ ([Wd]).  
\\
\\
3.2. {\em $L$-functions associated to the universal deformations.} Let $K$ be a knot in the $3$-sphere $S^3$ and  let $E_K$  denote the knot complement $S^3 \setminus K$. Let $G_K$ denote the knot group $\pi_1(E_K)$  of $K$, which has the following presentation:
$$
G_K=\left< g_1, \dots, g_n \; | \; r_1 = \cdots = r_{n-1} = 1 \right>. \leqno{(3.2.1)}
$$
Let $F$ be the free group on the words $g_1, \dots , g_n$ and let $\pi : \mathbb{Z}[F] \rightarrow \mathbb{Z}[G]$ be the natural homomorphism
of group rings. We write the same $g_i$ for the image of $g_i$ in $G_K$. 

Let $\overline{\rho} : G_K \rightarrow {\rm SL}_2(k)$ be an absolutely irreducible representation of $G_K$ over a perfect field $k$ with ${\rm char}(k) \neq 2$. Let ${\cal O}$ be a complete discrete valuation ring with residue field $k$ and let $\frak{CL}_{\cal O}$ be the category of complete local ${\cal O}$-algebras with residue field $k$. Let ${\boldsymbol \rho} : G_K \rightarrow {\rm SL}_2({\boldsymbol R}_{\overline{\rho}})$ be the universal deformation of $\overline{\rho}$ (Theorem 1.2.2).  We denote by the same ${\boldsymbol \rho}$ the induced algebra homomorphism $\mathbb{Z}[G_K] \rightarrow {\rm M}_2({\boldsymbol R}_{\overline{\rho}})$. Let $V_{\boldsymbol \rho}$ be the representation space (of column vectors) $({\boldsymbol R}_{\overline{\rho}})^{\oplus 2}$  of ${\boldsymbol \rho}$ on which $G_K$ acts from the left via ${\boldsymbol \rho}$.  We will compute the {\em twisted knot module}
$$H_*({\boldsymbol \rho}) := H_*(E_K; V_{\boldsymbol \rho})$$
 with coefficients in  $V_{\boldsymbol \rho}$ as the homology of the chain complex $C_*(W; V_{\boldsymbol \rho})$ of the CW complex $W$ attached to the presentation (3.2.1). The CW complex $W$ was given in Subsection 2.3.  Since $H_1(W;\mathbb{Z}) = H_1(E_K;\mathbb{Z}) = \langle t \rangle \simeq \mathbb{Z}$, we take  $\alpha : \pi_1(W) \rightarrow \langle t \rangle$  to be the abelianization map.

For a representation $\rho : G_K \rightarrow {\rm SL}_2(A)$, where $A$ is a Noetherian factorial domain, we define the {\em twisted Alexander invariant} $\Delta_K(\rho; t)$ of $K$ associated to $\rho$ by 
$$\Delta_K(\rho; t) := \Delta(W,\rho; t).$$
We note that $\Delta_K(\rho;t)$ coincides with the Reidemeister torsion of $E_K$ (or $W$) associated to the representation $\rho \otimes \alpha$ over $Q(A)(t)$ ([cf. FV; Proposition 2.2]). 
It is known  ([FV], [KL], [Ki]) that the relation between the twisted Alexander invariant $\Delta_K(\rho;t)$ and the initial Fitting ideals of $H_i(\rho \otimes \alpha) := H_i(E_K; \rho \otimes \alpha)$ $(i=0, 1)$  is given by \\
$$ \Delta_K(\rho; t)  \; \dot{=} \; \frac{\Delta_0(H_1(\rho \otimes \alpha))}{\Delta_0(H_0(\rho \otimes \alpha))}.  \leqno{(3.2.2)} $$
\\
Following Mazur's question 1 of [Ma3; page 440], we may ask the following \\
\\
{\bf Problem 3.2.3.}  Is $H_1({\boldsymbol \rho})$ a finitely generated and torsion ${\boldsymbol R}_{\overline{\rho}}$-module ? \\
\\
Here is our main theorem, which gives an affirmative answer to Problem 3.2.3 under some conditions using a twisted Alexander invariant of $K$.\\
\\
{\bf Theorem 3.2.4.} {\em Notations being as above, suppose that the following two conditions are satisfied}:\\
(1) {\em ${\boldsymbol R}_{\overline{\rho}}$ is a Noetherian integral domain.} \\
(2) {\em There is a deformation $\rho : G_K \rightarrow {\rm SL}_2(R)$ of $\overline{\rho}$, where $R \in \frak{CL}_{\cal O}$ is a Noetherian factorial domain, and $g \in G_K$ such that}\\
(2-1) $\det (\rho(g)-I) \not= 0$ {\em and} \\
(2-2) $\Delta_K(\rho;1) \neq 0$.
\\
{\em Then  $H_1({\boldsymbol \rho})$ is a finitely generated torsion ${\boldsymbol R}_{\overline{\rho}}$-module.  }
\\
\\
{\em Proof.}  We may assume that $g = g_n$ in the presentation (3.2.1) of $G_K$�D
We consider the following chain complex $ C_*({\boldsymbol \rho}) := C_*(W; V_{\boldsymbol \rho})$ ([Kw; 7.1]):
$$
0 \longrightarrow C_2({\boldsymbol \rho}) \overset{\partial_2}{\longrightarrow} C_1({\boldsymbol \rho}) \overset{\partial_1}{\longrightarrow} C_0({\boldsymbol \rho}) \longrightarrow 0,
$$
defined by
$$ \left\{ \begin{array}{l} C_0({\boldsymbol \rho}) := V_{\boldsymbol \rho},\\
C_1({\boldsymbol \rho}) := (V_{\boldsymbol \rho})^{\oplus n},\\
C_2({\boldsymbol \rho}) := (V_{\boldsymbol \rho})^{\oplus (n-1)},
\end{array} \right.
\;\;
\left\{ \begin{array}{l}
\partial_1 := ({\boldsymbol \rho}(g_1)-I,  \dots ,{\boldsymbol \rho}(g_n) -I ), \\
\partial_2 := ( \displaystyle{{\boldsymbol \rho} \circ \pi \left( \frac{\partial r_j}{\partial g_i} \right) }),
\end{array} \right.
$$
where $\frac{\partial}{\partial g_i} : \mathbb{Z}[F] \rightarrow \mathbb{Z}[F]$ denotes  the Fox derivative ([Fo]), and $\partial_2$ is regarded as a (big) $n \times (n-1)$ matrix whose $(i,j)$-entry is the $2 \times 2$ matrix ${\boldsymbol \rho} \circ \pi \left( \frac{\partial r_j}{\partial g_i} \right) $. 

 By the condition (2), let $\psi : {\boldsymbol R}_{\overline{\rho}} \rightarrow R$ be a morphism in $\frak{CL}_{\cal O}$ such that $\psi \circ {\boldsymbol \rho} \approx \rho$. Since 
$\psi(\det ({\boldsymbol \rho}(g_n)-I)) = \det(\rho(g_n) - I) \neq 0$ by (2-1),  we have $\det ({\boldsymbol \rho}(g_n)-I) \in Q({\boldsymbol R}_{\overline{\rho}})^{\times}$ by the condition (1)�DHence 
we have 
$$H_0({\boldsymbol \rho})\otimes_{{\boldsymbol R}_{\overline{\rho}}} Q({\boldsymbol R}_{\overline{\rho}}) =0.  \leqno{(3.2.4.1)} $$

Let $C_1'({\boldsymbol \rho})$ be the ${\boldsymbol R}_{\overline{\rho}}$-submodule of $C_1({\boldsymbol \rho})$ consisting of the first $(n-1)$ components so that $C_1({\boldsymbol \rho}) = C_1'({\boldsymbol \rho}) \oplus V_{\boldsymbol \rho}$ and let $\partial'_2$ be the $(n-1) \times (n-1)$ matrix obtained deleting the $n$-th row from $\partial_2$. Consider the ${\boldsymbol R}_{\overline{\rho}}$-homomorphism
$$
\partial'_2 : C_2({\boldsymbol \rho}) \longrightarrow C_1'({\boldsymbol \rho}).
$$
Then, by the definition (3.1.3) of the twisted Alexander invariant,  we have
$$
\Delta_K (\rho; 1) = \frac{ \psi(\det (\partial'_2))}{ \psi(\det ({\boldsymbol \rho}(g_n) - I)) }. \leqno{(3.2.4.2)}
$$
By the conditions (2-1),  (2-2) and  (3.2.4.2), we have $\det (\partial'_2) \in Q({\boldsymbol R}_{\overline{\rho}})^{\times}$. Hence we have
$$H_2({\boldsymbol \rho})\otimes_{{\boldsymbol R}_{\overline{\rho}}} Q({\boldsymbol R}_{\overline{\rho}}) =0.  \leqno{(3.2.4.3)} $$

 Since the Euler characteristic of $W$ is zero, we have
$$ \begin{array}{ll} 
\displaystyle{ \sum_{i=0}^3 (-1)^i \dim_{Q({\boldsymbol R}_{\overline{\rho}})} H_i({\boldsymbol \rho}) \otimes_{{\boldsymbol R}_{\overline{\rho}}}   Q({\boldsymbol R}_{\overline{\rho}})} & = \displaystyle{\sum_{i=0}^3 (-1)^i \dim_{Q({\boldsymbol R}_{\overline{\rho}})}  C_i({\boldsymbol \rho}) \otimes_{{\boldsymbol R}_{\overline{\rho}}}  Q({\boldsymbol R}_{\overline{\rho}})} \\
 & = \displaystyle{({\rm rank}_{{\boldsymbol R}_{\overline{\rho}}}\; V_{\boldsymbol \rho}) \sum_{i=0}^3 (-1)^i {\rm rank}_{\mathbb{Z}}\; C_i(W) }\\
 & = 0. 
 \end{array}
\leqno{(3.2.4.4)}
$$
Therefore, by (3.2.4.1), (3.2.4.3) and (3.2.4.4), we have 
$$ {\rm rank}_{{\boldsymbol R}_{\overline{\rho}}} \; H_1({\boldsymbol \rho})  = \dim_{Q({\boldsymbol R}_{\overline{\rho}})}  H_1({\boldsymbol \rho}) \otimes_{{\boldsymbol R}_{\overline{\rho}}}  Q({\boldsymbol R}_{\overline{\rho}}) = 0
$$
and hence $H_1({\boldsymbol \rho})$ is torsion over ${\boldsymbol R}_{\overline{\rho}}$. Since ${\boldsymbol R}_{\overline{\rho}}$ is Noetherian and $H_1({\boldsymbol \rho})$  is a quotient of a submodule of  $(V_{\boldsymbol \rho})^{\oplus n} = ({\boldsymbol R}_{\overline{\rho}})^{\oplus 2n}$, $H_1({\boldsymbol \rho})$ is Noetherian, in particular, finitely generated over ${\boldsymbol R}_{\overline{\rho}}$.  $\;\; \Box$\\
\\
It may be interesting to note that the condition (2-2) in Theorem 3.2.4 on a twisted Alexander polynomial is reminiscent of Kato's result in number theoretic situation ([Kt]), which asserts that the non-vanishing of the $L$-function at 1 of a modular form implies the finiteness of the Selmer module of the associated $p$-adic Galois representation.\\
As a special case of Theorem 3.2.4, the above proof  shows the following.\\
\\
{\bf Corollary 3.2.5.} {\em Notations being as above, suppose that the following two conditions are satisfied}:\\
(1) {\em ${\boldsymbol R}_{\overline{\rho}}$ is a Noetherian integral domain.} \\
(2) {\em There is $g \in G_K$ such that $\det (\overline{\rho}(g)-I) \not= 0$ and $\Delta_K(\overline{\rho};1) \neq 0$.}
\\
{\em Then  we have $H_1({\boldsymbol \rho}) = 0$. }\\
\\
{\it Proof.} By the assumptions, we have $\det({\boldsymbol \rho}(g_n)-I), \det(\partial_2') \in ({\boldsymbol R}_{\overline{\rho}})^{\times}$, from which we easily see that ${\rm Ker}(\partial_1) = {\rm Im}(\partial_2)$ and hence $H_1({\boldsymbol \rho}) = 0$. $\;\; \Box$\\
\\

Assume that ${\boldsymbol R}_{\overline{\rho}}$ is a Noetherian factorial domain and the condition (2) of Theorem 3.2.4.  When $H_1({\boldsymbol \rho})$ is a torsion ${\boldsymbol R}_{\overline{\rho}}$-module, we are interested in the invariant
$$ L_K({\boldsymbol \rho}) := \Delta_0(H_1({\boldsymbol \rho})),  \leqno{(3.2.6)}$$ 
which we call the {\em  $L$-function} of the knot $K$ associated to ${\boldsymbol \rho}$ (cf. Remark 3.2.8 (3) below). We note that it is a computable invariant by the following\\
\\
{\bf Proposition 3.2.7.} {\em Notations being as above, we have}
$$ L_K({\boldsymbol \rho}) \, \dot{=} \,  \Delta_2({\rm Coker}(\partial_2) ) $$  
\vspace{.2cm}
{\em Proof}.  This follows from the exact sequence of ${\boldsymbol R}_{\overline{\rho}}$-modules
$$ 0 \longrightarrow H_1({\boldsymbol \rho}) \longrightarrow {\rm Coker}(\partial_2) \longrightarrow V_{{\boldsymbol \rho}}=({\boldsymbol R}_{\overline{\rho}})^{\oplus 2} \longrightarrow 0$$
and Lemma 3.1.2 (2). $\;\; \Box$ \\
\\
{\bf Remark 3.2.8.} (1) The $L$-function $L_K({\boldsymbol \rho})$ is determined up to multiplication by a unit of ${\boldsymbol R}_{\overline{\rho}}$.\\
(2) When $H_*(\boldsymbol \rho) \otimes_{{\boldsymbol R}_{\overline{\rho}}} Q({\boldsymbol R}_{\overline{\rho}}) = 0$,  we have the Reidemeister torsion $\Delta_K({\boldsymbol \rho};1) \in Q({\boldsymbol R}_{\overline{\rho}})$ of $E_K$ associated to ${\boldsymbol \rho}$,  which is an invariant defined without indeterminacy. It may be non-trivial, even when $H_*(\boldsymbol \rho) = 0$.\\
(3) Our $L$-function $L_K({\boldsymbol \rho})$ may be seen as an analogue in knot theory of the algebraic $p$-adic $L$-function for the universal Galois deformation in number theory ([G]). 
In terms of [Ma3],  the ${\boldsymbol R}_{\overline{\rho}}$-module  $H_1({\boldsymbol \rho})$ gives a coherent torsion sheaf ${\cal H}_1({\boldsymbol \rho})$ on the universal deformation space ${\rm Spec}({\boldsymbol R}_{\overline{\rho}})$ and  $L_K({\boldsymbol \rho})$ gives a non-zero section of  ${\cal H}_1({\boldsymbol \rho})$. \\

We find the following necessary condition for the $L$-function $L_K({\boldsymbol \rho})$ to be non-trivial under a mild condition.\\
\\
{\bf Proposition 3.2.9.} {\em Assume that $\Delta_0(H_0({\boldsymbol \rho})) \; \dot{=} \;  1$. If $L_K({\boldsymbol \rho}) \; \dot{\neq} \; 1$, we have $\Delta_K(\overline{\rho};1) = 0$.}
\\
\\
{\em Proof.} By (3.2.2), (3.2.6) and our assumption, we have
$$ \Delta_K({\boldsymbol \rho};1) \; \dot{=} \;  L_K({\boldsymbol \rho}). \leqno{(3.2.9.1)}$$
Suppose  $L_K({\boldsymbol \rho}) \; \dot{\neq} \; 1$, which means $L_K({\boldsymbol \rho}) \in \frak{m}_{{\boldsymbol R}_{\overline{\rho}}}$. Let  $\varphi : {\boldsymbol R}_{\overline{\rho}} \rightarrow k$ be the homomorphism taking mod $\frak{m}_{{\boldsymbol R}_{\overline{\rho}}}$. Then, by the functorial property of the twisted Alexander invariant and (3.2.9.1) , we have
$$ \Delta_K({\overline{\rho}};1) = \varphi(\Delta_K({\boldsymbol \rho};1)) = \varphi(L_K({\boldsymbol \rho})) = 0. \;\; \Box $$
\\

Following Mazur's question 2 of [Ma3; page 440], we may ask the following\\
\\
{\bf Problem 3.2.10.} Investigate the order of the zeroes of $L_K({\boldsymbol \rho})$ on ${\rm Spec}({\boldsymbol R}_{\overline{\rho}})$ at prime divisors. \\
\\
In the next section, we verify Problem 3.2.10 affirmatively by some examples.\\
\\
{\bf Remark 3.2.11.} In [Ma3], Mazur works over a field $k$ (in fact, the field of complex numbers) and so the $L$-function discussed there is, in our terms, given by
$$ L_K({\boldsymbol \rho}_k) := \Delta_0(H_1({\boldsymbol \rho}_k)),$$
 where  ${\boldsymbol \rho}_k : G_K \rightarrow {\rm SL}_2({\boldsymbol R}_{\overline{\rho}} \otimes_{\cal O} k )$ is the representation  obtained by taking mod $\frak{m}_{\cal O}$ of  ${\boldsymbol \rho}$. Therefore our $L$-function $L_K({\boldsymbol \rho})$ in (3.2.6) is a finer object than $L_K({\boldsymbol \rho}_k)$.\\

\begin{center}
{\bf 4. Examples}
\end{center}

In this section, we discuss concrete examples of the universal deformations of some  representations of 2-bridge knot groups over finite fields and the associated $L$-functions. \\
 
 Let $K$ be a 2-bridge knot in the 3-sphere $S^3$, given as the Schubert form $B(m,n)$ where $m$ and $n$ are odd integers with $m>0, -m<n<m$ and ${\rm g.c.d}(m,n) =1$. The knot group $G_K$ is known to have a presentation of the form
$$ G_K = \langle g_1, g_2 \; | \; wg_1 = g_2w \rangle,$$
where $w$ is a word $w(g_1,g_2)$ of $g_1$ and $g_2$ which has the following symmetric form
$$ \begin{array}{l}
w = w(g_1,g_2) = g_1^{\epsilon_1} g_2^{\epsilon_2} \cdots g_1^{\epsilon_{m-2}} g_2^{\epsilon_{m-1}}, \\
\epsilon_i = (-1)^{[in/m]} = \epsilon_{m-i} \;\; ([\, \cdot \,] = \,\mbox{Gauss symbol}). \end{array}$$
We write the same $g_i$ for the image of the word $g_i$ in $G_K$.

Let $A$ be  a commutative ring with identity. For $a \in A^{\times}$ and $b \in A$, we consider two matrices $C(a)$ and $D(a, b)$ in ${\rm SL}_2(A)$ defined by
$$ C(a) := \left(  \begin{array}{cc} a & 1 \\ 0 & a^{-1} \end{array} \right), \;\;  D(a, b) := \left(  \begin{array}{cc} a & 0 \\  b & a^{-1} \end{array} \right)   $$
and we set
$$ W(a, b) := C(a)^{\epsilon_1} D(a, b)^{\epsilon_2}\cdots C(a)^{\epsilon_{m-2}} D(a, b)^{\epsilon_{m-1}}.$$
It is easy to see that there are (Laurent) polynomials $w_{ij}(t,u) \in \mathbb{Z}[t^{\pm},u]$ $(1\leq i,j \leq 2)$ such that $W(a,b) = (w_{ij}(a,b))$. Let $\varphi(t,u) := w_{11}(t,u) + (t^{-1} - t)w_{12}(t,u) \in \mathbb{Z}[t^{\pm},u]$. Then it is shown ([R]) that there is a unique polynomial $\Phi(x,u) \in \mathbb{Z}[x,u]$ such that
$$ \Phi(t+t^{-1}, u) = t^{l} \varphi(t,u)$$
for an integer $l$.

Let $k$ be a field with ${\rm char}(k) \neq 2$ and let ${\cal O}$ be a complete discrete valuation ring with residue field $k$. Let ${\cal X}(G_K)_k$ denote the character variety of $G_K$ over $k$. The proof of Proposition 1.4.1 of [CS] tells us that any ${\rm tr}(\sigma_{G_K}(g))$ ($g \in G_K$) is given as a polynomial of ${\rm tr}(\sigma_{G_K}(g_1)) (= {\rm tr}(\sigma_{G_K}(g_2)))$ and ${\rm tr}(\sigma_{G_K}(g_1g_2))$ with coefficients in $\mathbb{Z}$. In particular, the character algebra  ${\cal B}(G_K)_k$ is generated by ${\rm tr}(\sigma_{G_K}(g_1))$ and ${\rm tr}(\sigma_{G_K}(g_1g_2))$ over $k$. Let $x$ and $y$ denote the variables corresponding,  respectively,  to the coordinate functions ${\rm tr}(\sigma_{G_K}(g_1))$ and ${\rm tr}(\sigma_{G_K}(g_1g_2))$ on ${\cal X}(G_K)$. This variable $x$ is consistent with the variable $x$ of $\Phi(x,u)$ (and so causes no confusion). Since ${\rm tr}(C(a)) = a + a^{-1}$ and ${\rm tr}(C(a)D(a,b)) = a^2 + a^{-2} + b$, the coordinate variables $x$ and $y$ are related with $t$ and $u$  by 
$$x = t + t^{-1}, \;\; y = t^2 + t^{-2} + u = x^2 + u -2.$$
The following theorem is due to Le.\\
\\
{\bf Theorem 4.1} ([Le, Theorem 3.3.1]). {\em We have}
$$ {\cal X}(G_K)_k =  {\rm Spec} (k[x,y]/((y-x^2+2)\Phi(x, y-x^2+2))).$$
{\em Here, for a $k$-algebra $A$,  the $A$-rational points on $\Phi(x, y-x^2+2) = 0$ correspond bijectively to isomorphism classes  of absolutely irreducible representation $G_K \rightarrow {\rm SL}_2(A)$ except the finitely many intersection points with $y -x^2 +2 = 0$.}\\ 
\\ 
 {\bf Example 4.2.} (1) When $K$ is  the trefoil knot $B(3,1)$, we see $\Phi(x, y-x^2+2) = y-1$. \\
 (2) When $K$ is the figure eight knot $B(5,3)$,   we have $\Phi(x,y-x^2+2) = y^2 -(1+x^2)y + 2x^2 -1$. \\
 (3)  When $K := B(7,3)$, the knot $5_2$, we have  $\Phi(x,y-x^2+2) = y^3 - (x^2 + 1)y^2 + (3x^2 - 2)y - 2x^2 + 1$. \\
 \\
 By Theorem 4.1, we have the following \\
\\
{\bf Corollary 4.3.} {\em Let $\overline{\rho} : G_K \rightarrow {\rm SL}_2(k)$ be an absolutely irreducible representation so that $[\overline{\rho}]$ is a regular $k$-rational point of ${\cal X}(G_K)_k^{\rm a.i}$. Then we have}
$$ ({\cal B}(G_K)_k)_{[\overline{\rho}]}^{\wedge} \simeq k[[x - {\rm tr}(\overline{\rho}(g_1))]].$$
\\
So, by Theorem 2.2.4, we have\\
\\
{\bf Corollary 4.4.} {\em Let $\overline{\rho}$ be as in Corollary 4.3. Suppose that $\rho : G_K \rightarrow {\rm SL}_2({\cal O}[[x - \alpha]])$, where $\alpha$ is an element of ${\cal O}$ such that $\alpha \; \mbox{mod} \; \frak{m}_{\cal O} =   {\rm tr}(\overline{\rho}(g_1))$, is a deformation of $\overline{\rho}$ satisfying}
$$ {\rm tr}(\rho (g_1)) = x.   \leqno{(4.4.1)}$$
{\em Then the pair $({\cal O}[[x - \alpha]], \rho)$ is the universal deformation of $\overline{\rho}$.}
\\

In the following, we discuss some concrete examples, where $k$ will be a finite field $\mathbb{F}_p$ for some odd prime number $p$.\\
\\
{\em Convention}. Let $R$ be a complete local ring with residue field $R/\frak{m}_R = \mathbb{F}_p$. When the equation $X^2 = a$ for $a \in R$ has two simple roots in $R$, we denote by $\sqrt{a}$ for the "positive" solution, namely, $(\sqrt{a})^2 = a$ and $\sqrt{a} \; \mbox{mod}\; \frak{m}_R \in \{ 1, \dots , \frac{p-1}{2} \}$. \\
\\
{\bf Example 4.5.} (1) Let $K := B(3,1)$, the trefoil knot,  whose group is given by 
$$G_{K} = \langle g_1, g_2 \ | \ g_1 g_2 g_1 = g_2 g_1 g_2 \rangle. $$
We have ${\cal X}(G_K)_{\rm a.i}(k) = \{ (x,y) \in k^2 \; | \; y = 1\}$.\\

Let $k = \mathbb{F}_{3}$ and ${\cal O}= \mathbb{Z}_{3}$, and consider the following absolutely irreducible representation whose ${\rm PGL}_2(\mathbb{F}_3)$-conjugacy class corresponds to  
the regular $\mathbb{F}_3$-rational point  $(x, y) = (2,1)$ of  ${\cal X}(G_K)_{\rm a.i}$ (Proposition 2.1.3): 
$$\overline{\rho}_{1}: G_{K} \to  {\rm SL}_{2}(\mathbb{F}_{3}); \;\; 
\overline{\rho}_{1}(g_1) 
= \begin{pmatrix} 0 & 2 \\ 1 & 2 \end{pmatrix}, \ 
\overline{\rho}_{1}(g_2) 
= \begin{pmatrix} 2 & 2 \\ 1 & 0 \end{pmatrix}. $$
Let $\boldsymbol \rho_{1}: G_{K} \to {\rm SL}_{2}(\mathbb{Z}_{3}[[x - 2]])$ 
be the representation defined by 
$$ \begin{array}{c}
\boldsymbol \rho_{1}(g_1) = 
\begin{pmatrix} 
\frac{x + \sqrt{x^{2} - 3}}{2} & 
-1 \\ 
\frac{1}{4}  & 
\frac{x - \sqrt{x^{2} - 3}}{2} 
\end{pmatrix}, 
\\
\boldsymbol \rho_{1}(g_2) = 
\begin{pmatrix} 
\frac{x - \sqrt{x^{2} - 3}}{2} & 
-1 \\ 
\frac{1}{4}  & 
\frac{x + \sqrt{x^{2} - 3}}{2} 
\end{pmatrix}. 
\end{array}
$$
We see by the straightforward computation 
that $\boldsymbol \rho_{1}$ is indeed a  representation of $G_K$ and a deformation of $\overline{\rho}_{1}$ (see our convention).
 Moreover, we have ${\rm tr}(\boldsymbol \rho_{1}(g_1)) = x$, hence $\boldsymbol \rho_{1}$ satisfies the condition (4.4.1). 
Therefore $(\boldsymbol R_{\overline{\rho}_{1}} = \mathbb{Z}_{3}[[x - 2]], \boldsymbol \rho_{1})$ is the universal 
deformation of $\overline{\rho}_{1}$. 

We easily see that $\Delta_0(H_0({\boldsymbol \rho}_1)) \; \dot{=}\; 1$ and  $\Delta_K(\overline{\rho}_1;t) = 1 + t^2$, hence, $\Delta_K(\overline{\rho}_1;1) = 2 \neq 0$.
Therefore, by Proposition 3.2.9,  we have 
$$H_{1}(\boldsymbol \rho_{1}) = 0, \ L_K(\boldsymbol \rho_{1}) \; \dot{=}\;  1. $$
(2)  Let $K := B(5,3)$, the figure eight knot,  whose group is given by 
$$G_{K} = \langle g_1, g_2 \ | \ g_1 g_2^{-1} g_1^{-1} g_2 g_1 = g_2 g_1 g_2^{-1} g_1^{-1} g_2 \rangle. $$
We have ${\cal X}(G_K)_{\rm a.i}(k) = \{ (x,y) \in k^2 \; | \; y^2 - (1+x^2)y + 2x^2 -1 = 0 \} \setminus \{(\pm \sqrt{5},3)\}$.\\

Let $k = \mathbb{F}_{7}$ and ${\cal O} = \mathbb{Z}_{7}$, and consider the following absolutely irreducible representation
whose ${\rm PGL}_2(\mathbb{F}_7)$-conjugacy  class corresponds to the regular $\mathbb{F}_7$-rational points $(x,y) = (5,5)$ of ${\cal X}(G_K)_{\rm a.i}$: 
$$ \begin{array}{l}
\overline{\rho}_{2}: G_{K} \to {\rm SL}_{2}(\mathbb{F}_{7}); \; 
\overline{\rho}_{2}(g_1) = 
\begin{pmatrix} 0 & 6 \\ 1 & 5 \end{pmatrix}, 
\overline{\rho}_{2}(g_2) = 
\begin{pmatrix} 5 & 6 \\ 1 & 0 \end{pmatrix}. 
\end{array}
$$
Let $\boldsymbol \rho_{2}: G_{K} \to {\rm SL}_{2}(\mathbb{Z}_{7}[[x + 2]])$ 
be the representation defined by 
$$ \begin{array}{c}
\boldsymbol \rho_{2}(g_1) = 
\begin{pmatrix} 
 \frac{x + \sqrt{\frac{x^{2} - 5 + u(x)}{2}}}{2} & 
 -1 \\ 
 -\frac{x^{2} - 3 - u(x)}{8} & 
 \frac{x - \sqrt{\frac{x^{2} - 5 + u(x)}{2}}}{2} 
\end{pmatrix}, 
\\
\boldsymbol \rho_{2}(g_2) = 
\begin{pmatrix} 
 \frac{x - \sqrt{\frac{x^{2} - 5 + u(x)}{2}}}{2} & 
 -1 \\ 
 -\frac{x^{2} - 3 - u(x)}{8} & 
 \frac{x + \sqrt{\frac{x^{2} - 5 + u(x)}{2}}}{2} 
\end{pmatrix}, 
\end{array}
$$ 
where $u(x) := \sqrt{(x^{2} - 1) (x^{2} - 5)}$. 
We see by the straightforward computation
that $\boldsymbol \rho_{2}$ is indeed a  representation of $G_K$ and a deformation of $\overline{\rho}_{2}$.
Moreover, we have ${\rm tr}(\boldsymbol \rho_{2}(g_1)) = x$, 
hence $\boldsymbol \rho_{2}$ satisfies the condition (4.4.1). 
Therefore $(\boldsymbol R_{\overline{\rho}_{2}} = \mathbb{Z}_{7}[[x + 2]], \boldsymbol \rho_{2})$ 
is the universal deformation of $\overline{\rho}_{2}$. 

We easily see that $\det(\overline{\rho}_{2}(g_2) - I) = 4 \neq 0$ and that $\Delta_{K}(\overline{\rho}_{2}; t) = t^{-2} + 4t^{-1} + 1$, hence, $\Delta_{K}(\overline{\rho}_{2}; 1) = 6 \neq 0$.
 Therefore, by Corollary 3.2.5, we have 
$$H_{1}(\boldsymbol \rho_{2}) = 0, \; L_K(\boldsymbol \rho_{2}) \; \dot{=}\;  1. $$
(3)  Let $K := B(7,3)$,   the knot $5_{2}$, whose group is given by 
$$G_{K} = \langle g_1, g_2 \ | \ g_1 g_2 g_1^{-1} g_2^{-1} g_1 g_2 g_1 = g_2 g_1 g_2 g_1^{-1} g_2^{-1} g_1 g_2 \rangle. $$
We have ${\cal X}(G_K)_{\rm a.i}(k) = \{ (x,y) \in k^2 \; | \; y^3 - (x^2 + 1)y^2 +(3x^2 -2)y - 2x^2 + 1 = 0 \} \setminus \{ (\pm \sqrt{\frac{7}{2}}, \frac{3}{2})\}$.\\

Firstly, let $k = \mathbb{F}_{11}$ and ${\cal O} = \mathbb{Z}_{11}$, and consider the following absolutely irreducible representation whose ${\rm PGL}_2(\mathbb{F}_{11})$-conjugacy class corresponds
to the regular $\mathbb{F}_{11}$-rational point $(x,y) = (5,5)$ of ${\cal X}(G_K)_{\rm a.i}$: 
$$\overline{\rho}_{3}: G_{K} \to {\rm SL}_{2}(\mathbb{F}_{11}); \ 
\overline{\rho}_{3}(g_1) 
= \begin{pmatrix} 5 & 10 \\ 1 & 0 \end{pmatrix}, \ 
\overline{\rho}_{3}(g_2) 
= \begin{pmatrix} 5 & 1 \\ 10 & 0 \end{pmatrix}. $$
Let $\alpha := \frac{3-\sqrt{5}}{2}, \xi := \frac{4- \sqrt{5}}{4} \in \mathbb{Z}_{11}$ so that $\alpha \; {\rm mod}\; 11 = 5, \xi \; \mbox{mod}\; 11 = 0 \in \mathbb{F}_{11}$. 
Let $s = s(x)$ be the unique solution in $\mathbb{Z}_{11}[[x - \alpha]]$ satisfying the equation
$$64s^{3} - 16(2x^{2} + 5)s^{2} + 4(x^{4} + 9x^{2} + 6)s - 4x^{4} - 6x^{2} - 1 = 0 \leqno{(4.5.1)} $$ 
and
$$ s(\alpha) = \xi. \leqno{(4.5.2)}$$
Such an $s(x)$ is proved, by Hensel's lemma ([Se2; II, $\S 4$, Proposition 7]), to exist uniquely. Now, let $\boldsymbol \rho_{3}: G_{K} \to {\rm SL}_{2}(\mathbb{Z}_{11}[[x - \alpha]])$ 
be the representation defined by 
$$ \begin{array}{c}
\boldsymbol \rho_{3}(g_1) = 
\begin{pmatrix} 
\frac{x + \sqrt{x^{2} - 4s(x)}}{2} & 
-1 \\ 
-s(x) + 1  & 
\frac{x - \sqrt{x^{2} - 4s(x)}}{2} 
\end{pmatrix}, 
\\
\boldsymbol \rho_{3}(g_2) = 
\begin{pmatrix} 
\frac{x + \sqrt{x^{2} - 4s(x)}}{2} & 
1 \\ 
s(x) - 1 & 
\frac{x - \sqrt{x^{2} - 4s(x)}}{2}  
\end{pmatrix}.
\end{array}
$$
 We can verify by (4.5.1) that $\boldsymbol \rho_{3}$ is indeed a representation of $G_K$ and by (4.5.2) that $\boldsymbol \rho_{3}$ is a deformation of $\overline{\rho}_{3}$. 
Moreover, we have  ${\rm tr}(\boldsymbol \rho_{3}(g_1)) = x$, 
hence $\boldsymbol \rho_{3}$ satisfies the condition (4.4.1). 
Therefore $(\boldsymbol R_{\overline{\rho}_{3}} = \mathbb{Z}_{11}[[x - \alpha]], \boldsymbol \rho_{3})$ is the universal 
deformation of $\overline{\rho}_{3}$. 

Consider the 11-adic lifting $\rho_{3}: G_{K} \to {\rm SL}_2(\mathbb{Z}_{11})$ of $\overline{\rho}_{3}$ 
defined  by ${\boldsymbol \rho}_3|_{x= 5}$:
$$
\rho_{3}(g_1) = 
\begin{pmatrix} 
 \frac{5 + \sqrt{25 - 4\mu}}{2} & 
 -1 \\ 
 -\mu + 1 & 
 \frac{5 - \sqrt{25 - 4\mu}}{2} 
\end{pmatrix}, 
$$
$$
\rho_{3}(g_2) =
\begin{pmatrix} 
 \frac{5 + \sqrt{25 - 4\mu}}{2} & 
 1\\
 \mu - 1 & 
 \frac{5 - \sqrt{25 - 4\mu }}{2} 
\end{pmatrix},
$$
where $\mu$ is the unique solution in $\mathbb{Z}_{11} $ satisfying (4.5.1) with $x=5$ and $\mu \; \mbox{mod}\; 11 = 0$. Then we easily see that 
$\det(\rho_{3}(g_2) - I) = -3 \neq 0$, 
and that
$\Delta_{K}(\rho_{3}; t) = -2 \{-8\mu^2 + 58\mu - 52 + 5t + (-8\mu^2 + 58\mu -52)t^2 \}$, hence,  
$\Delta_{K}(\rho_{3}; 1) = -2(-16\mu^2 + 116 \mu -99) \neq 0$. Therefore, by Theorem 3.2.4, $H_1(\boldsymbol \rho_3)$ is a finitely generated torsion $\mathbb{Z}_{11}[[x - \alpha]]$-module. 

 We let $r := g_1 g_2 g_1^{-1} g_2^{-1} g_1 g_2 g_1g_2^{-1}g_1^{-1}g_2g_1g_2^{-1}g_1^{-1}g_2^{-1}$ and set 
 $$ 
\partial_2 = ({\boldsymbol \rho}_3( \frac{\partial r}{\partial g_1}), {\boldsymbol \rho}_3( \frac{\partial r}{\partial g_2}))  = ({\boldsymbol a}_1, {\boldsymbol a}_2, {\boldsymbol a}_3, {\boldsymbol a}_4).  
$$
By the computer calculation, we find that all 2-minors of $\partial_2$ are given by
$$ \begin{array}{l}
\det ({\boldsymbol a}_1, {\boldsymbol a}_2) =   2(x-2)\{ 4(s-1)x^2 + x -4(2s-1)^2\},\\
\det ({\boldsymbol a}_1, {\boldsymbol a}_3) =  - \frac{1}{2} \{ 4(s-1)x^4 - 2(8s^2 - 2s -5)x^2 + 4(s-1)x \\
      \;\;\;\;\;\;\;\;\;\;\;\;\;\;\;\;\;\;\;\;  + (4s-3)(12s-5)\}(x-2-\sqrt{x^2 -4s}),\\
\det ({\boldsymbol a}_1, {\boldsymbol a}_4) =  4(s-1)x^4 - 8(s-1)x^3 - 4(4s^2 - 5s + 2)x^2 \\
 \;\;\;\;\;\;\;\;\;\;\;\;\;\;\;\;\;\;\;\;  + 4(8s^2 - 7s + 2)x -(4s-1)^2,\\
\det ({\boldsymbol a}_2, {\boldsymbol a}_3) =  -\{  4(s-1)x^4 - 8(s-1)x^3 - 4(4s^2 - 5s + 2)x^2 \\
 \;\;\;\;\;\;\;\;\;\;\;\;\;\;\;\;\;\;\;\;   + 4(8s^2 - 7s + 2)x -(4s-1)^2 \},\\
\det ({\boldsymbol a}_2, {\boldsymbol a}_4)  =  2\{ 4(s-1)x^2 + x -4(2s-1)^2\}(x-2 + \sqrt{x^2 - 4s}),\\
\det ({\boldsymbol a}_3, {\boldsymbol a}_4)  =  2(x-2)\{ 4(s-1)x^2 + x -4(2s-1)^2\}.
\end{array}
\leqno{(4.5.3)}
$$
By (4.5.3) and the computer calculation,  we find that $x = \alpha$ ($s(\alpha) = \frac{4 - \sqrt{5}}{4}$) gives a common zero of all 2-minors of $\partial_2$ and and their derivatives and is not a common zero of the third order derivatives of all 2-minors. Hence the greatest common divisor of all 2-minors  is $(x-\alpha)^2$. Therefore, by Proposition 3.2.7, we have 
$$H_{1}(\boldsymbol \rho_{3}) \simeq \mathbb{Z}_{11}\oplus \mathbb{Z}_{11}, \; L_K(\boldsymbol \rho_{3})\; \dot{=}\;  (x - \alpha)^2. $$

\vspace{.5cm}

Secondly, let $k = \mathbb{F}_{19}$ and ${\cal O} = \mathbb{Z}_{19}$, and consider the following absolutely irreducible representation whose ${\rm PGL}_2(\mathbb{F}_{19})$-conjugacy class corresponds to
the regular $\mathbb{F}_{19}$-rational point $(x,y) = (6,6)$ of ${\cal X}(G_K)_{\rm a.i}$: 
$$\overline{\rho}_{4}: G_{K} \to {\rm SL}_{2}(\mathbb{F}_{19}); \ 
\overline{\rho}_{4}(g_1) 
= \begin{pmatrix} 14 & 1 \\ 1 & 11 \end{pmatrix}, \ 
\overline{\rho}_{4}(g_2) 
= \begin{pmatrix} 11 & 1 \\ 1 & 14 \end{pmatrix}. $$
Let $\beta := \frac{3+\sqrt{5}}{2}, \zeta := \frac{7+\sqrt{5}}{8} \in \mathbb{Z}_{19}$ so that $\beta \; \mbox{mod}\; 19 = 6, \zeta \; \mbox{mod}\; 19 = 2 \in \mathbb{F}_{19}$. 
Let  $v = v(x)$ be the unique solution in $\mathbb{Z}_{19}[[x - \beta]]$ satisfying the equation 
$$64v^{3} - 16(x^{2} + 7)v^{2} + 28(x^2 + 2)v  - 12x^2 - 7 = 0 \leqno{(4.5.4)} $$ 
and
$$ v(\beta) = \zeta. \leqno{(4.5.5)}$$
Such a $v(x)$ is proved, by Hensel's lemma ([Se2; II, $\S 4$, Proposition 7]), to exist uniquely.
Now, let $\boldsymbol \rho_{4}: G_{K} \to {\rm SL}_{2}(\mathbb{Z}_{19}[[x - \beta]])$ 
be the representation defined by 
$$ \begin{array}{c}
\boldsymbol \rho_{4}(g_1) = 
\begin{pmatrix} 
\frac{x + \sqrt{x^{2} - 4v(x)}}{2} & 
1 \\ 
v(x) - 1  & 
\frac{x - \sqrt{x^{2} - 4v(x)}}{2} 
\end{pmatrix}, 
\\
\boldsymbol \rho_{4}(g_2) = 
\begin{pmatrix} 
\frac{x - \sqrt{x^{2} - 4v(x)}}{2} & 
1 \\ 
v(x) - 1 & 
\frac{x + \sqrt{x^{2} - 4v(x)}}{2}  
\end{pmatrix}.
\end{array}
$$
 We can verify by (4.5.4) 
that $\boldsymbol \rho_{4}$ is indeed a representation of $G_K$ and by (4.5.5) that $\boldsymbol \rho_{4}$ is a deformation of $\overline{\rho}_{4}$.
 Moreover, we have ${\rm tr}(\boldsymbol \rho_{4}(g_1)) = x$, 
hence $\boldsymbol \rho_{4}$ satisfies the condition (4.4.1). 
Therefore $(\boldsymbol R_{\overline{\rho}_{4}} = \mathbb{Z}_{19}[[x - \beta]], \boldsymbol \rho_{4})$ is the universal 
deformation of $\overline{\rho}_{4}$. 

Consider the 19-adic lifting $\rho_{4}: G_{K} \to {\rm SL}_2(\mathbb{Z}_{19})$ of $\overline{\rho}_{4}$ 
defined  by ${\boldsymbol \rho}_4|_{x= 6}$:
$$
\rho_{4}(g_1) = 
\begin{pmatrix} 
 \frac{6 + \sqrt{36 - 4\nu}}{2} & 
 1 \\ 
 \nu - 1 & 
 \frac{6 - \sqrt{36 - 4\nu}}{2} 
\end{pmatrix}, 
$$
$$
\rho_{4}(g_2) =
\begin{pmatrix} 
 \frac{6 - \sqrt{36 - 4\nu}}{2} & 
 1\\
 \nu - 1 & 
 \frac{6 + \sqrt{36 - 4\nu }}{2} 
\end{pmatrix},
$$
where $\nu$ is the unique solution in  $\mathbb{Z}_{19} $ satisfying  (4.5.4) with $x=6$ and  $\nu \; \mbox{mod}\; 19 = 2$. Then we easily see that 
$\det(\rho_{4}(g_2) - I) = -4 \neq 0$, 
and that
$\Delta_{K}(\rho_{4}; t) = -2 \{-8\nu^2 + 80\nu - 74 + 6t + (-8\nu^2 + 80\nu -74)t^2 \}$, hence,  
$\Delta_{K}(\rho_{4}; 1) = -2(-16\nu^2 + 160 \nu - 142) \neq 0$. Therefore, by Theorem 3.2.4, $H_1(\boldsymbol \rho_4)$ is a finitely generated torsion $\mathbb{Z}_{19}[[x - \beta]]$-module. 

We  set 
 $$ 
\partial_2 = ({\boldsymbol \rho}_4( \frac{\partial r}{\partial g_1}), {\boldsymbol \rho}_4( \frac{\partial r}{\partial g_2}))  = ({\boldsymbol b}_1, {\boldsymbol b}_2, {\boldsymbol b}_3, {\boldsymbol b}_4).  
$$
By the computer calculation, we find that all 2-minors of $\partial_2$ are given by
$$ \begin{array}{l}
\det ({\boldsymbol b}_1, {\boldsymbol b}_2) =   2(x-2)\{ 4(v-1)x^2 + x -4(2v-1)^2\},\\
\det ({\boldsymbol b}_1, {\boldsymbol b}_3) =  - \frac{1}{2} \{ 4(v-1)x^2 - 4(v-1)x - (4v-3)^2\}\sqrt{x^2 - 4v} \\
\det ({\boldsymbol b}_1, {\boldsymbol b}_4) =  4(v-1)x^4 - (8v-9)x^3 - 2(8v^2 - 10v + 5)x^2 \\
 \;\;\;\;\;\;\;\;\;\;\;\;\;\;\;\;\;\;\;\;  + 4(8v^2 - 9v + 3)x -(4v-3)^2\\
  \;\;\;\;\;\;\;\;\;\;\;\;\;\;\;\;\;\;\;\;     - (x-2)\{ 4(v-1)x^2 + x -4(2v-1)^2\}\sqrt{x^2 - 4v},\\
\det ({\boldsymbol b}_2, {\boldsymbol b}_3) =  -\{  4(v-1)x^4 - (8v-9)x^3 - 2(8v^2 - 10v + 5)x^2 \\
 \;\;\;\;\;\;\;\;\;\;\;\;\;\;\;\;\;\;\;\;   + 4(8v^2 - 9v + 3)x -(4v-3)^2 \}\\
   \;\;\;\;\;\;\;\;\;\;\;\;\;\;\;\;\;\;\;\;       - (x-2)\{ 4(v-1)x^2 + x -4(2v-1)^2\}\sqrt{x^2 - 4v},\\
\det ({\boldsymbol b}_2, {\boldsymbol b}_4)  =  2\{ 4(v-1)x^2 + x -4(2v-1)^2\}\sqrt{x^2 - 4v},\\
\det ({\boldsymbol b}_3, {\boldsymbol b}_4)  =  2(x-2)\{ 4(v-1)x^2 + x -4(2v-1)^2\}.
\end{array}
\leqno{(4.5.6)}
$$
By (4.5.6) and the computer calculation,  we find that $x = \beta$ ($v(\beta) = \frac{7+\sqrt{5}}{8}$) is a common zero of all 2-minors of $\partial_2$ and their derivatives and is not a common zero of the third order derivatives of all 2-minors. Hence the greatest common divisor of all 2-minors  is $(x-\beta)^2$. Therefore, by Proposition 3.2.7, we have 
$$H_{1}(\boldsymbol \rho_{4}) \simeq \mathbb{Z}_{19}\oplus \mathbb{Z}_{19}, \; L_K(\boldsymbol \rho_{4})\; \dot{=}\;  (x - \beta)^2. $$
\\

We see that all examples above answer Problems 3.2.3 affirmatively and answer Problem 3.2.10 concretely.\\

\begin{flushleft}
{\bf References}\\
{[Ca]} H. Carayol, Formes modulaires et repr\'{e}sentations galoisiennes \`{a} valeurs dans un anneau local complet, $p$-adic monodromy and the Birch and Swinnerton-Dyer conjecture (Boston, MA, 1991), Amer. Math. Soc.,  Contemp. Math., {\bf 165}, 1994, 213--237. \\
{[CS]} M. Culler, P. Shalen, Varieties of group representations and splittings of $3$-manifolds, Ann. of Math. {\bf 117}  (1983), 109--146.\\
{[Fo]} R. Fox, Free differential calculus. I. Derivation in the free group ring, Ann. of Math. (2)  {\bf 57}  (1953). 547--560.\\
{[FV]} S. Friedl, S. Vidussi, A survey of twisted Alexander polynomials, The mathematics of knots, 45--94,
Contrib. Math. Comput. Sci., 1, Springer, Heidelberg, 2011.\\
{[Fu]} T. Fukaya, Hasse zeta functions of non-commutative rings, J. Algebra  {\bf 208}  (1998),  no. 1, 304--342. \\
{[G]} R. Greenberg,  Iwasawa theory and $p$-adic deformations of motives, Motives (Seattle, WA, 1991),  Proc. Sympos. Pure Math., {\bf 55}, Part 2, Amer. Math. Soc., Providence, RI, 1994, 193--223.\\
{[Ha]} S. Harada, Modular Representations of Fundamental Groups and Associated Weil-type Zeta Functions, Thesis, Kyushu University, 2008.\\
{[H1]} H. Hida, Galois representations into ${\rm GL}_2(\mathbb{Z}_p[[X]])$ attached to ordinary cusp forms, Invent. Math.  {\bf 85},  (1986),  no. 3, 545--613. \\
{[H2]} H. Hida, Iwasawa modules attached to congruences of cusp forms, Ann. Sci. \'{E}cole Norm. Sup. (4)  {\bf 19},  (1986),  no. 2, 231--273.\\
{[Hi]} J. Hillman, Algebraic invariants of links, Second edition. Series on Knots and Everything, {\bf 52}, World Scientific Publishing Co. Pte. Ltd., Hackensack, NJ, 2012.\\
{[I]} K. Iwasawa, On $\mathbb{Z}_l$-extensions of algebraic number fields, Ann. of Math. (2)  {\bf 98} (1973), 246--326. \\
{[Kt]}  K. Kato, $p$-adic Hodge theory and values of zeta functions of modular forms,  Cohomologies $p$-adiques et applications arithm\'{e}tiques. III. Ast\'{e}risque  No. {\bf 295},  (2004),  117--290.
\\
{[Kw]} A. Kawauchi, A survey of knot theory, Translated and revised from the 1990 Japanese original by the author. Birkh\"{a}user Verlag, Basel, 1996.\\
{[Ki]} T. Kitano, Twisted Alexander polynomial and Reidemeister torsion,
Pacific J. Math. {\bf 174}  (1996), no. 2, 431--442. \\
{[KL]} P. Kirk, C. Livingston, Twisted Alexander invariants, Reidemeister torsion, and Casson-Gordon invariants, 
Topology  {\bf 38}  (1999),  no. 3, 635--661. \\
{[Ko]} K. Kodama, Knot program,  available at http://www.math.kobe-u.ac.jp/~kodama/knot.html \\
{[Ku]} M. Kurihara, Iwasawa theory and Fitting ideals, J. Reine Angew. Math.  {\bf 561}  (2003), 39--86.\\
{[Le]} T. Le, Varieties of representations and their subvarieties of cohomology jumps for certain knot
groups,  Russian Acad. Sci. Sb. Math. {\bf 78}, (1994), 187--209.\\
{[LM]} A. Lubotzky, A. Magid, Varieties of representations of finitely generated groups, Mem. Amer. Math. Soc. vol. {\bf 58}, no. {\bf 336}, 1985. \\
{[Ma1]} B. Mazur, Remarks on the Alexander Polynomial, available at  http://www.math.harvard.edu/~mazur/older.html \\
{[Ma2]} B. Mazur,  Deforming Galois representations, Galois groups over $\mathbb{Q}$,  Math. Sci. Res. Inst. Publ., {\bf 16}, Springer, (1989), 385--437.\\
{[Ma3]} B. Mazur, The theme of $p$-adic variation,  Mathematics: frontiers and perspectives, AMS, RI, (2000), 433--459. \\
{[MW]} B. Mazur, A. Wiles, Class fields of abelian extensions of $\mathbb{Q}$,  Invent. Math.  {\bf 76}  (1984),  no. 2, 179--330. \\
{[Mo]} M. Morishita, Knots and Primes - An introduction to Arithmetic Topology, Universitext, Springer, London, 2012.\\
{[MTTU]} M. Morishita, Y. Takakura, Y. Terashima, J. Ueki, On the universal deformations for ${\rm SL}_2$-representations of knot groups, to appear in Tohoku Math. J.\\
 {[Na]} K. Nakamoto, Representation varieties and character varieties, Publ. Res. Inst. Math. Sci.  {\bf 36}  (2000),  no. 2, 159--189.\\
{[Ny]} L. Nyssen, Pseudo-repr\'{e}sentations, Math. Ann., {\bf 306}  (1996), no.2, 257--283.\\

{[O1]} T. Ochiai, Control theorem for Greenberg's Selmer groups of Galois deformations, J. Number Theory  {\bf 88}  (2001),  no. 1, 59--85. \\

{[O2]} T. Ochiai, On the two-variable Iwasawa main conjecture, Compos. Math.  {\bf 142},  (2006),  no. 5, 1157--1200.\\

{[Po]} J. Porti, Torsion de Reidemeister pour les vari\'{e}t\'{e}s hyperboliques,  Mem. Amer. Math. Soc.  {\bf 128}  (1997),  no. 612.\\

{[Pr]}  C. Procesi, The invariant theory of $n \times n$ matrices,  Advances in Math.  {\bf 19} (1976), no. 3, 306--381.\\

{[PS]} J. Przytycki, A. Sikora, On skein algebras and $Sl_2(\mathbb{C})$-character varieties, Topology, {\bf 39} (2000), 115--148.\\

{[R]} R. Riley, Nonabelian representations of 2-bridge knot groups, Quar. J. Oxford, {\bf 35} (1984), 191--208.\\

{[Sa]} K. Saito, Character variety of representations of a finitely generated group in ${\rm SL}_2$, Topology and Teichm\"{u}ller spaces (Katinkulta, 1995),  World Sci. Publ., River Edge, NJ, (1996), 253--264.\\

{[Se1]} J.-P. Serre, Classes des corps cyclotomiques (d'apres K. Iwasawa), S\'{e}minaire Bourbaki, Vol. 5,  Exp. No. {\bf 174}, 11 pages, Soc. Math. France, Paris, 1958.\\

{[Se2]} J.-P. Serre, Corps locaux, Publications de l'Universit\'{e} de Nancago, No. VIII. Hermann, Paris, 1968. \\

{[Se3]} J.-P. Serre, Cohomologie Galoisienne, Lecture Notes in Mathematics, Vol. {\bf 5}. Springer-Verlag, Berlin-New York, 1973. \\

{[Ta]} R. Taylor, Galois representations associated to Siegel modular forms of low weight, Duke Math. J. {\bf 63} (1991), no. 2, 281--332. \\

{[Ti]} J. Tilouine, Deformations of Galois representations and Hecke algebras, Published for The Mehta Research Institute of Mathematics and 
Mathematical Physics, Allahabad; by Narosa Publishing House, New Delhi, 1996. \\

{[Wd]} M. Wada, Twisted Alexander polynomial for finitely presentable groups, Topology  {\bf 33}  (1994),  no. 2, 241--256. \\

{[Ws]}L. Washington, Introduction to cyclotomic fields, Second edition. Graduate Texts in Mathematics, {\bf 83}, Springer-Verlag, New York, 1997. \\
\end{flushleft}

\noindent
{\small  Takahiro Kitayama:\\
Graduate School of Mathematical Sciences,\\
the University of Tokyo,\\
3-8-1 Komaba, Meguro-ku, Tokyo 153-8914, JAPAN\\
e-mail:kitayama@ms.u-tokyo.ac.jp\\
\\
Masanori Morishita:\\
Faculty of Mathematics, \\
Kyushu University \\
744, Motooka, Nishi-ku, Fukuoka, 819-0395, JAPAN \\
e-mail: morisita@math.kyushu-u.ac.jp\\
\\
Ryoto Tange:\\
Faculty of Mathematics, \\
Kyushu University \\
744, Motooka, Nishi-ku, Fukuoka, 819-0395, JAPAN \\
e-mail: r-tange@math.kyushu-u.ac.jp\\
\\
Yuji Terashima:\\
Department of Mathematical and Computing Sciences, \\
Tokyo Institute of Technology, \\
2-12-1 Oh-okayama, Meguro-ku, Tokyo 152-8551, JAPAN\\
e-mail: tera@is.titech.ac.jp}

\end{document}